\documentclass[journal]{IEEEtran_local}
\usepackage[top=0.4in, bottom=0.2in, left=0.47in, right=0.15in]{geometry}
\usepackage{booktabs}
\usepackage{multirow}
\usepackage{tabularx}
\usepackage{geometry}
\usepackage{enumitem}
\usepackage{bbm}
\usepackage{mathtools}

\usepackage{subfigure}

\usepackage{cite}

\usepackage{amssymb}

\usepackage[pdftex]{graphicx}

\usepackage{amsmath}
\usepackage{amsfonts}
\usepackage{amsthm}

\theoremstyle{plain}

\newtheorem{prop}{Proposition}

\newtheorem{remark}{Remark}

\usepackage{algorithm}
\usepackage{algorithmic}
\usepackage{multirow}  
\usepackage{booktabs}  

\usepackage{array}

\usepackage{graphicx} 
\ifCLASSOPTIONcompsoc
\usepackage[caption=false,font=normalsize,labelfont=sf,textfont=sf]{subfig}
\else
\usepackage[caption=false,font=footnotesize]{subfig}
\fi
\usepackage{color}

\usepackage{stfloats}
\usepackage{bm}

\renewcommand{\arraystretch}{1.3}

\usepackage{nomencl}
\usepackage{ifthen}
\renewcommand{\nomgroup}[1]{%
	\ifthenelse{\equal{#1}{A}}{\item[\emph{\textbf{Set and Index}}]}{%
		\ifthenelse{\equal{#1}{B}}{\item[\emph{\textbf{Parameters}}]}{%
			\ifthenelse{\equal{#1}{C}}{\item[\emph{\textbf{Variables}}]}
		}
	}
}
\makenomenclature

\makeatletter
\renewcommand{\subsection}{\@startsection{subsection}{2}{\z@}%
{0.5ex}
{0.2ex}
	{\normalfont\normalsize\itshape}} 
\makeatother

\makeatletter
\renewcommand{\section}{\@startsection{section}{1}{\z@}%
	{1.0ex plus 0.5ex minus 0.5ex}%
	{0.4ex plus 0.2ex}%
	{\normalfont\normalsize\centering\scshape}}
\makeatother

\begin{document}
	\allowdisplaybreaks[4]
%
\title{Distributionally Robust Planning of Hydrogen-Electrical Microgrids for Sea Islands}
%
%
%

\author{Yuchen Dong, Zhengsong Lu, Xiaoyu Cao, Zhengwen He, Tanveer Hossain Bhuiyan, Bo Zeng
\thanks{Yuchen Dong and Zhengwen He are with the School of Management and the Key Lab of the Ministry of Education for Process Management \& Efficiency Engineering, Xi’an Jiaotong University, Xi’an, Shaanxi 710049 China (e-mail: dongyuchen@stu.xjtu.edu.cn; zhengwenhe@mail.xjtu.edu.cn).}
\thanks{Zhengsong Lu and Bo Zeng are with the Department of Industrial Engineering, University of Pittsburgh, Pittsburgh, Pennsylvania, 15261 USA (e-mail: zs.lu@pitt.edu; bzeng@pitt.edu).}
\thanks{Xiaoyu Cao is with the School of Automation Science and Engineering and the Ministry of Education Key Lab for Intelligent Networks and Network Security, Xi’an Jiaotong University, Xi’an, Shaanxi 710049 China (e-mail: cxykeven2019@xjtu.edu.cn).}
\thanks{Tanveer Hossain Bhuiyan is with the Department of Mechanical, Aerospace, and Industrial Engineering, The University of Texas at San Antonio, San Antonio, Texas, 78249, USA (e-mail: Tanveer.Bhuiyan@utsa.edu).}
}%
\maketitle

\begin{abstract}
This paper presents a distributionally robust planning  method for hydrogen-electrical microgrids over islands, where the cross-island energy exchange is supported by a maritime hydrogen transport network.  
This planning problem is complicated due to heterogeneous off-shore wind-driven uncertainties (i.e., renewable power, transport availability, demand fluctuations, and grid faulting), a subset of which exhibit endogenous uncertainty, as they can be affected by proactive measures (e.g., grid hardening) or infrastructure investment. To capture these features, a two-stage distributionally robust optimization (DRO) model is developed considering decision-dependent uncertainty (DDU), which encompasses variation of the underlying distributional ambiguity due to the change of the  first stage decisions. Notably, the complete recourse property is missing, which is often neglected in existing DRO studies. Nevertheless, different from the case for land-based microgrids,  this issue is critical and fundamental for sea island systems due to their particular physical and logistical requirements.  To address these issues, we develop a C\&CG algorithm that is customized with strong cutting planes to handle DRO with a varying DDU ambiguity set and feasibility requirements.  Numerical results demonstrate the cost-effectiveness and resilience of the proposed planning framework, along with the nontrivial improvements of the algorithm in both solution accuracy and computational efficiency.
\end{abstract}
\vspace{-5pt}

%
\IEEEpeerreviewmaketitle

\printnomenclature
\section{Introduction}
%
%
%
%

\IEEEPARstart{S}{ea} islands, as ecologically and economically vital regions, play an increasingly prominent role in global sustainable development, with growing demands for clean and reliable power supply  to  advance sustainable tourism, low-carbon resource development, and resilient island communities. 
However, due to their geographic isolation and lack of interconnection with mainland power grids, these islands typically rely on stand-alone distribution systems. 
Traditional power supply strategies, depending on the transport of fossil fuels, often result in high operational costs, poor reliability, and environmental  deterioration. To address these challenges, it is essential to shift toward renewable-based microgrids that offer greater self-sufficiency, cost viability, and environmental friendliness.

Renewable energy microgrids on sea islands have received increasing attentions in recent years,
which can be categorized into methodological developments \cite{kuznia2013stochastic,li2022noncooperative,zhou2019capacity,li2018optimal, sadiq2024towards,teng2025distributed} and practical implementations \cite{de2010isolated,silva2019optimal}.
From a methodological perspective,  \cite{kuznia2013stochastic} presented a two-stage stochastic programming (SP) model to support system planning considering random renewable energy sources (RESs).  
\cite{li2022noncooperative} proposed a noncooperative game-theoretic framework for multi-operator planning in islanded systems, focusing on renewable integration under imperfect information. \cite{zhou2019capacity} developed a storage control and sizing strategy to reduce reliance on costly diesel generation. To achieve reliable and economical operations,  \cite{li2018optimal} formulated a chance-constrained microgrid scheduling model accounting for reserve uncertainty. 
On the practical side, the viability of renewable-based islanded microgrids has been demonstrated by hybrid microgrid deployments in Len{\c{c}}{\'o}is island \cite{de2010isolated} and the Azores archipelago \cite{silva2019optimal}, both of which are grid-isolated yet rich in natural energy resources. 

Despite these promising demonstrations, the development of sea-island microgrids remains fundamentally constrained by the spatial, ecological, and  infrastructural  characteristics that differentiate them from stand-alone systems.  
Specifically, many inhabited islands (referred to as load islands \cite{sui2020day}) face  strong regulations on large-scale renewable equipment due to land-use restrictions and landscape preservation.   
Given these considerations,  nearby resource islands can provide more favorable conditions for  renewable generation, co-located with associated industrial activities. This spatial separation between power generation and consumption naturally necessitates a  multi-microgrid configuration, through which stand-alone cross-island energy exchange is supported by diverse delivery pathways. 

Among various technical options for cross-island energy exchange, mobile energy delivery using dedicated carriers has been considered, especially when  submarine cable deployment is impractical due to high capital cost and installation complexity \cite{yang2023hybrid}.  
For example,  the use of electric vessels has been proposed by \cite{sui2020day} and \cite{mahmud2019real} for mobile electricity delivery via onboard batteries. 
Nevertheless, given  limitations in energy density, scalability, and long-duration storage performance, battery-based solutions are inadequate for energy transmission among multi-island power systems. On the other hand, it has been recognized that hydrogen presents a more effective solution under such a situation \cite{wang2023optimization,sui2023optimal}.   
With higher gravimetric energy density, hydrogen can be stored long-term in compressed or liquefied form with minimal loss. 
Beyond storage, it also serves as a cross-sectoral energy vector that links renewable power with storage and fuel-cell (FC) systems, while enabling temporal decoupling between energy production and consumption.
\cite{sui2023optimal} proposed a variable-speed scheduling strategy for hydrogen carrier vessels to realize cross-island energy transmission. While such operational strategies address hydrogen dispatch in the short term, the system-level planning and configuration of hydrogen-electrical (HE) multi-microgrids over sea islands has yet to be studied.

Such planning decisions are particularly challenging due to the inherent uncertainty in long-term infrastructure deployment over sea islands. 
These uncertainties include RES's variability, demand fluctuations, vessel availability and network disruptions.
A prominent source of these uncertainties is offshore wind,  exhibiting nonlinear effects on both energy supply and system resilience. For example, power generation is feasible only within turbine-specific cut-in and cut-out thresholds, while extreme wind events may suspend output and increase the risk of equipment failure or $N-k$ cascading disruptions.
Notably, a subset of these uncertainties exhibit decision-dependency (referred as DDU), meaning that their behavior is directly affected by  planning decisions. For instance, deploying generation infrastructure on resource islands induces local electricity demand from operational equipment, control systems, and on-site personnel.
Similarly, decisions on grid reinforcement and protection schemes influence the vulnerability of distribution lines to extreme weather, thereby shaping the occurrence probability and severity of faulting events \cite{li2023distributionally,pianco2025decision}. These characteristics necessitate a comprehensive  modeling framework that explicitly represents the fundamental connection between infrastructure decisions and uncertainty realizations.

To effectively hedge the complex uncertainties in sea island systems, a robust planning framework is essential,   
as the geographic isolation, limited monitoring infrastructure, and high environmental variability lead to scarce or unreliable data or information.
In this context, distributionally robust optimization (DRO) has been adopted as a powerful approach within a two-stage decision-making framework, where irreversible investment decisions are made in the first stage and operations follow in the second. 
Note that DRO seeks solutions that perform well in the expected sense across a family of plausible distributions, which are collectively referred to as an ambiguity set. 
For instance, \cite{xie2022sizing} proposed a DRO model with a Wasserstein-based ambiguity set to co-optimize renewable generation and energy storage sizing in isolated microgrids. 
Similarly, \cite{zhou2023distributionally} developed a moment-based DRO model for microgrid energy management under renewable uncertainty, which was reformulated as a tractable mixed-integer second-order cone program (SOCP). 
Actually, DRO’s advantages become especially salient under DDU, which is highly relevant for sea island systems. 
Line reinforcement decisions can shape the ambiguity sets used to model failure probability distributions, thereby affecting network resilience.
Moreover, demand shifts on resource islands, induced by the placement of renewable energy infrastructure, influence both the distributional parameters and the sample space.

While DRO offers an effective modeling framework, its computational tractability remains a significant challenge,  particularly under DDU context. There are two types of computational strategies to solve two-stage DRO exactly. \cite{li2023distributionally} and \cite{mohajerin2018data,zhao2018data,hanasusanto2018conic} solved its equivalent reformulation as a single-level (nonlinear) mixed integer program, which may encounter significant computational challenges in practice due to a large number of variables and constraints introduced. Another strategy decomposes the DRO model using the master-subproblem framework (e.g., \cite{bansal2018decomposition,gamboa2021decomposition, qu2025distributionally,yang2025robust,duque2022distributionally,pianco2025decision}) that often have a better solution capacity.  
We mention that a recent primal-based decomposition algorithm, which explicitly computes the worst-case distribution and the associated performance, has been developed in \cite{lu2024two} that achieves a superior scalability and efficiency over existing ones.  
Additionally, most solution frameworks rely on a complete recourse assumption, under which all first-stage decisions are assumed to yield feasible recourse actions for any realization of uncertainty. It is worth pointed  out that this assumption  is fundamentally unrealistic in sea island microgrids, due to limited access to external diesel supply and constrained storage capacity resulting from geographic remoteness.

To address the aforementioned planning challenges associated with HE multi-microgrids over sea islands,  this paper leverages  strong two-stage DRO modeling and computational tools and customizes them considering particular features and enhancements.  Specifically, infrastructure investment is determined in the first stage to minimize capital cost, while the  operational energy dispatch  are optimized in the second stage to minimize the worst-case expected operating and maintenance costs. One feature is a maritime hydrogen transport network that enables cross-island energy exchange. Another one is the consideration of both environmental and operational uncertainties and their impact on the feasibility issue, including RES variability, load fluctuations, hydrogen transport availability, and network contingencies. Most of them lack precise probabilistic characterizations due to data scarcity, and some of them are clearly DDU factors. In comparison to the current literature, this paper makes the following contributions:

\begin{enumerate}[left=0em]
	\item A comprehensive planning model for HE sea-island systems is first presented to capture the spatial-temporal decoupling of power supply and demand, facilitated by hydrogen acting as a cross-island energy carrier. In contrast to prior studies that focus on dispatch or rely on simplified transport, we further incorporate  a cycle-based multi-trip transport network that is jointly planned with energy infrastructure.
	\item We discretize nonlinear off-shore wind-related uncertainty into multiple levels to  capture its impact on renewable generation availability, hydrogen vessel operation, and $N$–$k$ network contingency. 
	Building upon this, we  propose a DDU formulation both on ambiguity set (e.g., line faults mitigated by reinforcement decisions) and sample space (e.g., induced loads from generation siting on resource islands). Most existing studies simplify key decision-dependent structures or consider only limited uncertainty dimensions to retain tractability.
	\item The developed two-stage DRO framework 
	accounts for an often fundamental requirement underlying sea island systems that load shedding either is not allowed or must be bounded under any possible scenario.  
	Actually, although such feasibility issues happen frequently in real life, they have been largely ignored or only crudely approximated in existing DRO literature. 

	\item    
    A strong primal-based decomposition algorithm, i.e.,  a variant of column-and-constraint generation (C\&CG) to handle DDU-DRO (referred to as C\&CG-DRO), has been developed with tailored enhancements. It produces effective cutting-planes to tackle the challenges arising from the absence of the aforementioned complete recourse assumption, and achieves extraordinary computational performance. Indeed, it often drastically outperforms the state-of-the-art methods by more than one order of magnitude.
\end{enumerate}
The remainder of this paper is organized as follows. Section \ref{SM} formulates the microgrids planning model. The two-stage DDU-DRO model is proposed in Section \ref{TS_DDUDRO}. Section \ref{PDDRO} develops a primal-based decomposition algorithm. Section \ref{NS} presents numerical experiments. The conclusion and discuss are presented in Section \ref{CONC}. 

\section{System Modeling}
\label{SM}
\setlength{\jot}{1pt} 
\setlength{\abovedisplayskip}{5pt}
\setlength{\belowdisplayskip}{4pt}
\setlength{\abovedisplayshortskip}{2pt}
\setlength{\belowdisplayshortskip}{2pt}
\begin{figure}[!t]
	\centering
	\includegraphics[width=4in]{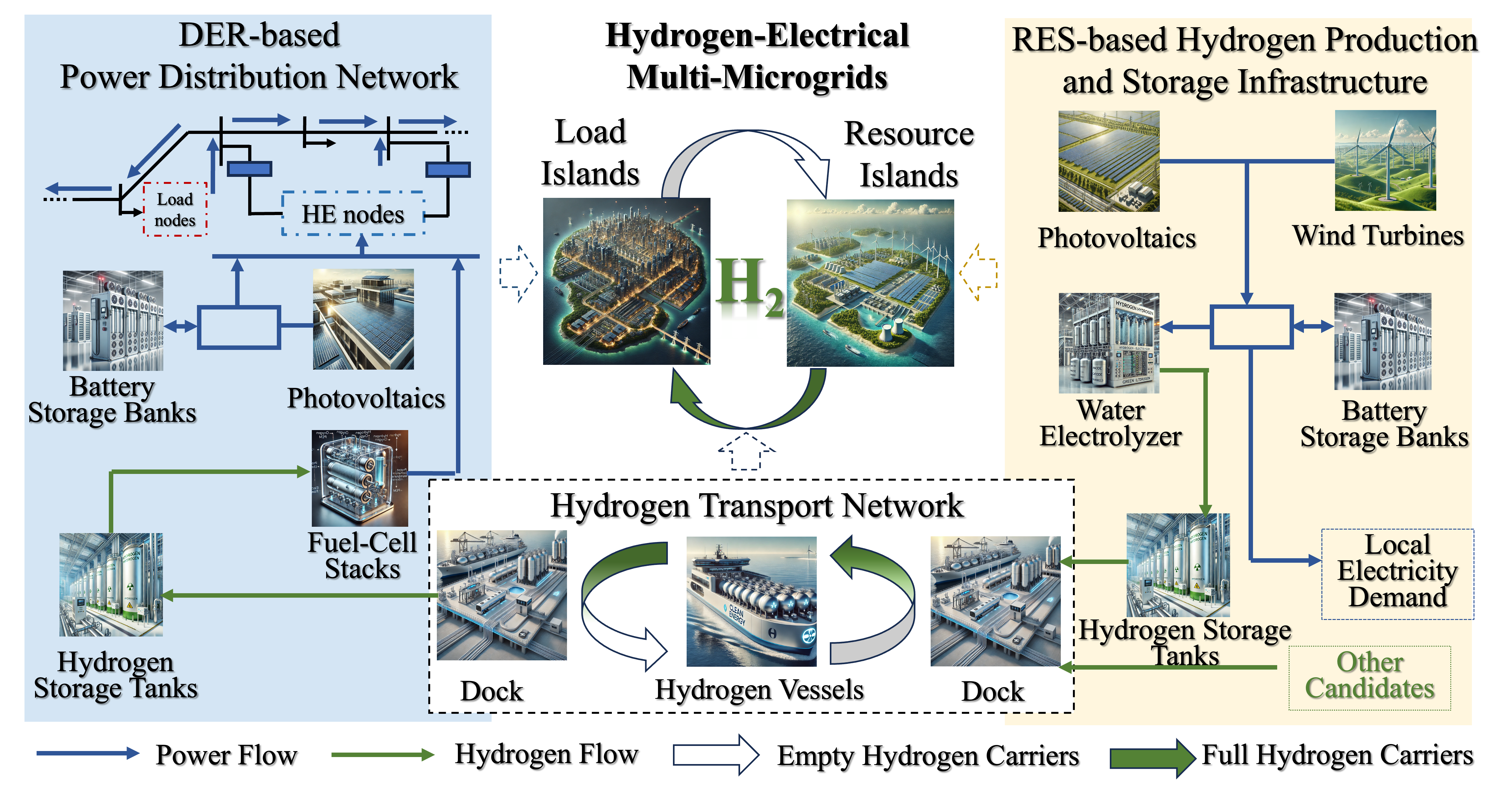}
	\caption{System architecture of hydrogen-electrical  sea-island microgrids with cross-island hydrogen transport and  energy conversion.} 
	\label{fig:system} 
\end{figure}
Fig. \ref{fig:system}   illustrates the proposed hydrogen-electrical (HE) sea-island microgrids, which is tailored for sea islands lacking feasible subsea cable connections. The system supports their decarbonization  by enabling deep renewable integration under spatial and infrastructural constraints.
In this architecture, resource islands focus on renewable hydrogen production, while load islands are designated for hydrogen utilization and electricity supply.
As the energy vector linking generation and consumption across islands, hydrogen is transported through a dedicated shipping network comprising storage facilities and hydrogen vessels that realize cross-island delivery. 

\subsection{Resource Islands: Renewable Hydrogen Production} 
Leveraging local renewable energy, the resource islands  (denoted by $s \in \mathcal{S}_r$) serve as decentralized hubs for  hydrogen production.  
 Within each resource island, multiple potential installation sites (indexed by $j\in\mathcal{N}_s$, for $s\in\mathcal{S}_r$) can be equipped with RES-based hydrogen production and storage (HP\&S)  infrastructure (i.e., electrolyzers, hydrogen storage tanks (HSTs), and battery storage banks (BSBs)). 
 The RES technologies, including wind turbines (WT) and photovoltaics (PV), are denoted by $o \in \{wt, pv\}$ with corresponding capacity decisions $\overline{P}_{o,j}$. 
For technologies in HP\&S units, their capacity planning decisions  are denoted by  $\overline{P}_{\mathfrak{d},j}, \ \mathfrak{d} \in \{\mathrm{elz, hst, bsb}\}$, respectively. 

\subsubsection{Operational Constraints for Power-to-Hydrogen (P2H)  Conversion} 
Note that locally generated renewable electricity is converted into hydrogen via electrolyzers, with surplus energy charged in BSBs and discharged when needed to maintain continuous operation. 
The P2H system is designed to achieve energy self-sufficiency at each site, whose operations are described as follows: 
 \setlength{\arraycolsep}{-0.2em}
\begin{eqnarray}
	&&0 \leq p_{o,j}^t \leq \tilde{\delta}_{o,j}^t\overline{P}_{o,j}, \ \forall j\in\mathcal{N}_{s}, o\in\{ wt \!,pv\}, \forall t\in\mathcal{T} \label{RE1} \\
	&&p_{\mathrm{elz},j}^t \leq \min\{\overline{P}_{\mathrm{elz},j},\textstyle\sum\nolimits_{o\in \{wt,pv\}}p_{o,j}^t + p_{\mathrm{dis},j}^t - p_{\mathrm{ch},j}^t - p^t_{\mathrm{l},j}\},\nonumber\\
	&& \qquad\qquad\qquad\qquad\qquad\qquad\qquad\forall j\in\mathcal{N}_{s},\forall t\in\mathcal{T} \label{RE2}\\
	&&m_{\mathrm{elz},j}^t=\eta_{\mathrm{p-h}}\chi_{\mathrm{p-h}} p_{\mathrm{elz},j}^t, \ \forall j\in\mathcal{N}_{s},\forall t\in\mathcal{T} \label{RE3}
\end{eqnarray}
where $p_{o,j}^t$ represents the power generated from WT or PV. The RES-based electricity output varies due to natural fluctuations, as captured in \eqref{RE1} by the resource availability factor $\tilde{\delta}_{o,j}^t$. The electricity supplied to the electrolyzers ($p_{\mathrm{elz},j}^t$) is constrained by both the electrolyzer capacity ($\overline{P}_{\mathrm{elz},j}$) and the available electricity ) from RESs and BSBs ($ p_{\mathrm{dis},j}^t - p_{\mathrm{ch},j}^t$) after meeting local demand ($p^t_{\mathrm{l},j})$, as formulated in \eqref{RE2}. Hydrogen production $m_{\mathrm{elz},j}^t$, described in \eqref{RE3}, is determined by the electricity input $p_{\mathrm{elz},j}^t$, the electrolysis efficiency $\eta_{\mathrm{p-h}}$, and the  conversion coefficient $\chi_{\mathrm{p-h}}$. 

\subsubsection{Dynamic Operations for Energy Storage Systems (ESSs)}
The HSTs perform the function as  buffers between on-site hydrogen production and periodic off-island transportation, mitigating production fluctuations and ensuring a stable supply for scheduled deliveries. 
Their operational constraints are described as follows: 
\begin{eqnarray}	
	&&(1\!-\!\varpi_{\mathrm{hst}})\overline{P}_{\mathrm{hst},j}\leq G_{\mathrm{hst},j}^{t+1}\! = \!(1\!-\!\rho_{\mathrm{hst}})G_{\mathrm{hst},j}^{t}\!+\!\Delta_t(\eta_{\mathrm{hst+}}m_{\mathrm{elz},j}^t\nonumber\\ 
	&&-{m_{j,t}^{\mathrm{dout}}}/{\eta_{\mathrm{hst-}}} )\leq \overline{P}_{\mathrm{hst},j},\ \forall j\in\mathcal{N}_{s},\forall t\in\mathcal{T} \label{RH1}
\end{eqnarray}
where $G_{\mathrm{hst},j}^t$ represents the hydrogen inventory level, whose  dynamical  balance is modeled in \eqref{RH1}, accounting for leakage losses ($\rho_{\mathrm{hst}}$), hydrogen inflow from electrolysis $m_{\mathrm{elz},j}^t$, and hydrogen extraction $m_{j,t}^{\mathrm{dout}}$ for shipping transportation, with filling/releasing efficiencies ($\eta_{\mathrm{hst+}}$/$\eta_{\mathrm{hst-}}$). Also, safe storage operations of HSTs are  bounded by the minimum/maximum storage capacity [$(1-\varpi_{\mathrm{hst}})\overline{P}_{\mathrm{hst},j}$,$\overline{P}_{\mathrm{hst},j}$]. 

Moreover, the BSBs charge during RES generation surplus and discharge to meet electrolyzer hydrogen demand, stabilizing HP\&S systems, as governed by the following constraints: 
\begin{eqnarray}	
	&&(1-\varpi_{\mathrm{bsb}})\overline{P}_{\mathrm{bsb},j}\leq  E_{\mathrm{bsb},j}^{t+1}=(1-\rho_{\mathrm{bsb}})E_{\mathrm{bsb},j}^{t} + \nonumber\\ 
	&&\!\Delta_t\left(\eta_{\mathrm{bsb+}}p_{\mathrm{ch},j}^t - {p_{\mathrm{dis},j}^t}/{\eta_{\mathrm{bsb-}}}\right) \leq \overline{P}_{\mathrm{bsb},j} \ \forall j\in\mathcal{N}_{s},\forall t\in\mathcal{T}\label{RBSB1}\\
	&&0\leq p_{\mathrm{\mathfrak{e}},j}^t\leq \varrho^{\mathrm{bsb}}_{\mathfrak{e}}  \overline{P}_{\mathrm{bsb},j},\ \forall \mathfrak{e} \in\{\mathrm{ch},\mathrm{dis}\}, \forall j\in\mathcal{N}_{s},\forall t\in\mathcal{T}   \label{RBSB3} \\
&&\textstyle \sum_{t\in\mathcal{T}}\left(\eta_{\mathrm{bsb+}}p_{\mathrm{ch},j}^t - {p_{\mathrm{dis},j}^t}/{\eta_{\mathrm{bsb-}}}\right)\!\Delta_t = 0, \forall j\in\mathcal{N}_s  \label{RBSB5}
\end{eqnarray}
where $E_{\mathrm{bsb},j}^{t}$ represents the BSBs' state of charge (SOC), whose balance is modeled in \eqref{RBSB1}, considering leakage losses ($\rho_{\mathrm{bsb}}$), charging/discharging efficiency ($\eta_{\mathrm{bsb+}}$/$\eta_\mathrm{{bsb-}}$). Also, $E_{\mathrm{bsb},j}^t$ should be with its operational range [$(1-\varpi_{\mathrm{bsb}})\overline{P}_{\mathrm{bsb},j}$,$\overline{P}_{\mathrm{bsb},j}$]. 
The charging and discharging power limits in \eqref{RBSB3} are defined by the coefficients ($\varrho^{\mathrm{bsb}}_{\mathrm{ch}}$ and $\varrho^{\mathrm{bsb}}_{\mathrm{dis}}$). As in  \eqref{RBSB5}, the end-period SOC equals its initial value.
\subsection{Load Islands: Hydrogen Utilization and Power Supply}
Load islands (represented by  $ s \in \mathcal{S}_d$) are geographically isolated distribution systems that face  high electricity demand but insufficient local energy supply. 
The limited availability of suitable land and regulatory constraints further restrict large-scale renewable deployment. 
To address this challenge, a hybrid strategy combines RES generation with FC-based power output from  imported hydrogen, while ESSs balance fluctuations and maximize the utilization of intermittent renewables to ensure stable electricity supply. 
\subsubsection{FC-based Hydrogen-to-Power (H2P) System}
FCs serve as a key component in H2P conversion, ensuring reliable power supply when PV generation is insufficient. The active/reactive power output of FCs is constrained by the following equations: 
\begin{eqnarray}	
	&&0 \leq p^t_{\mathrm{fc},j} = \eta_{\mathrm{h-p}}  \chi_{\mathrm{h-p}} m^t_{\mathrm{fc},j} \leq \overline{P}_{\mathrm{fc},j},\ \forall j\in\mathcal{N}_{s},\forall t\in\mathcal{T}  \label{LFC1}\\
	&&\tan \underline{\vartheta}_{j} p^t_{\mathrm{fc},j} \leq q^t_{\mathrm{fc},j} \leq \tan \overline{\vartheta}_{j} p^t_{fc,j}, \forall j\in\mathcal{N}_{s},\forall t\in\mathcal{T} \label{LFC4}
\end{eqnarray}
where  \eqref{LFC1} establishes the relationship between the active power output ($p^t_{\mathrm{fc},j}$) of FCs, limited by its capacity $\overline{P}_{\mathrm{fc},j}$,  and the hydrogen supply ($m^t_{\mathrm{fc},j}$) from the HST, where $ \eta_{\mathrm{h-p}} $ and $ \chi_{\mathrm{h-p}} $ are conversion efficiency parameters. \eqref{LFC4} further constrains the reactive power output ($p^t_{\mathrm{fc},j}$) within a range determined by the  power factor angles [$\underline{\vartheta}_{j},\overline{\vartheta}_{j}$], respectively. 
Since the operation of FCs relies on a stable hydrogen supply, the following constraints govern hydrogen storage dynamics:  
\begin{eqnarray}	
	&& G_{\mathrm{hst},j}^{t+1}=\Delta_t\left(\eta_{\mathrm{hst+}}m_{j,t}^{\mathrm{din}}- {m_{\mathrm{fc},j}^t}/{\eta_{\mathrm{hst-}}}\right) +(1-\rho_{\mathrm{hst}})G_{\mathrm{hst},j}^{t}, \nonumber\\ 
	&&\qquad\qquad\qquad\qquad\qquad\qquad\qquad \forall j\in\mathcal{N}_{s},\forall t\in\mathcal{T}  \label{LG1}
\end{eqnarray}
where storage balance constraints \eqref{LG1} consider hydrogen inflow from external supply ($m_{j,t}^{\mathrm{din}}$) to support its consumption from FCs ($m_{\mathrm{fc},j}^t$). 

\subsubsection{Power Distribution Network} 
Let $\mathcal{N}^a_s$ denote the set of nodes that constitute the grid of load island $s$.
A series of linearized  \textit{Distflow} constraints  (i.e., \eqref{PDN1}-\eqref{PDN6}) govern the flow of active/reactive power ($fp_{i,j}^t/fq_{i,j}^t$), as well as the nodal voltage levels ($ U_j^t$). 
\begin{eqnarray}
	&&\textstyle\sum\nolimits_{h\in\Phi(j)}fp_{j,h}^t-\sum\nolimits_{i\in\Upsilon(j)} fp_{i,j}^t=p_j^t, \ \forall j\in \mathcal{N}^a_s, \forall t\in\mathcal{T}\label{PDN1}\\
	&&\textstyle\sum\nolimits_{h\in\Phi(j)}fq_{j,h}^t-\sum\nolimits_{i\in\Upsilon(j)} fq_{i,j}^t=q_j^t, \ \forall j\in \mathcal{N}^a_s, \forall t\in\mathcal{T}\label{PDN22}\\
		&&U_i^t-U_j^t=({r_{i,j}fp_{i,j}^t+x_{i,j}fq_{i,j}^t})/{U_0}, \ \forall (i,j)\in \mathcal{E}, \forall t\in\mathcal{T} \  \label{PDN3} \\
		&&-\overline{FP}_{i,j}\leq fp_{i,j}^t\leq \overline{FP}_{i,j}, \ \forall (i,j)\in \mathcal{E}, \forall t\in\mathcal{T} \label{PDN4}\\
		&&-\overline{FQ}_{i,j}\leq fq_{i,j}^t\leq \overline{FQ}_{i,j}, \ \forall (i,j)\in \mathcal{E}, \forall t\in\mathcal{T} \label{PDN5}\\
	&&\underline{U}_j\leq U_j^t\leq \overline{U}_j, \ \forall j\in\mathcal{N}^a_s, \forall t\in\mathcal{T} \label{PDN6}	\\
	&&p_j^t=
	\left\{
	\begin{array}{lll}
		\ -{p}_{\mathrm{l},j}^t, \ \forall j\in \mathcal{N}^a_s\backslash\mathcal{N}_{s}, \forall t\in\mathcal{T}\\
		\ p_{pv,j}^t +p_{\mathrm{fc},j}^t + p_{\mathrm{dis},j}^t - p_{\mathrm{ch},j}^t-p_{\mathrm{l},j}^t,  \forall j\in \mathcal{N}_{s}, \forall t\in\mathcal{T}\\
	\end{array}
	\right.  \label{PDN7}\\
	&&q_j=
	\left\{
	\begin{array}{lll}
		\ -{q}_{\mathrm{l},j}^t, \ \forall j\in \mathcal{N}^a_s\backslash\mathcal{N}_{s}\\
		\ q_{\mathrm{fc},j}^t-q_{\mathrm{l},j}^t,  \ \forall j\in \mathcal{N}_{s}, \forall t\in\mathcal{T}
	\end{array}
	\right. \label{PDN8}
\end{eqnarray}
where the nodal active power $p_j^t$, as in \eqref{PDN7}, is defined based on the type of node. For non-source nodes ($j\in\mathcal{N}^a_s\setminus\mathcal{N}_s$),  active power is solely determined by local load ${p}_{\mathrm{l}, j}^t$. For source nodes ($j\in\mathcal{N}_s$),  active power accounts for RES generation, FCs output, battery charging/discharging, and local load. Similarly, as for   the reactive power in \eqref{PDN8}, the non-source nodes account for the reactive load ${q}_{\mathrm{l}, j}^t$, while the source nodes account for the reactive power compensation. 

\subsection{Hydrogen Transport Network: Cross-Island Energy Carrier}
To effectively integrate the energy planning and operations of an island cluster (denoted by set $\mathcal{S} = \mathcal{S}_d\cup\mathcal{S}_r$), we establish a structured maritime hydrogen transport network (HSN) that functions as a cross-island energy carrier, ensuring efficient and reliable hydrogen transportation. Within this network, hydrogen vessels (denoted by $v\in\mathcal{V}$) serve as dynamic energy links, continuously bridging the hydrogen supply-demand gap across islands. 
\subsubsection{Cycle-Based Multi-Trip Planning} 
To optimize transportation planning, we adopt a repetitive cycle-based approach, where the entire planning horizon is structured into multiple identical cycles. Each cycle represents a transportation period in which vessel operations follow a multi-trip routing and  scheduling framework. In each cycle, key binary decision variables include vessel allocation ($b_{v,s}$), $k$-th trip activation ($z_{v}^k$)  and route selection ($y_{ss',v}^k$), subject to the below constraints:
 \setlength{\arraycolsep}{-0.3em}
 \begin{eqnarray}
	&& \textstyle\sum\nolimits_{s\in\mathcal{S}_{d}}b_{v,s}\leq a^{\mathrm{ves}}_v,\ \forall v \in\mathcal{V} \label{HTR1}\\
	&&  z_{v}^k \leq a^{\mathrm{ves}}_v, \ k \in \mathcal{K},\  \forall v \in \mathcal{V}  \label{HTR2_1}\\
	&&\textstyle \sum\nolimits_{v\in\mathcal{V}} b_{v, s}\leq \overline{N}_{s}^{\mathrm{ves}}, \ \forall s \in\mathcal{S}_d \label{HTR2}\\
	&& z_{v}^{k+1} \leq z_{ v}^{k}, \ k \in \mathcal{K}, \forall v \in \mathcal{V}  \label{HTR3}\\
	&& y_{s,s',v}^k \leq z_{v}^k, \ \forall (s,s') \in\mathcal{S}\times\mathcal{S}, \forall k \in \mathcal{K}, \forall v \in \mathcal{V}\label{HTR4}\\
	&&\textstyle \sum\nolimits_{s'\in\mathcal{S}}y_{s,s', v}^k \leq b_{v,s} +  \sum\nolimits_{s''\in\mathcal{S}}y_{s'',s,v}^k,  \forall v \!\in\!\mathcal{V}, \forall k \!\in\! \mathcal{K}, \forall s \!\in\!\mathcal{S}_d \label{HTR5}\\
	&& \textstyle\sum\nolimits_{s'\in\mathcal{S}} y_{s,s',v}^k \! =\! \sum\nolimits_{s''\in\mathcal{S}}y_{s'',s,v}^k, \forall s\! \in\!\mathcal{S}_r, \forall k \!\in\! \mathcal{K}, \forall v \!\in\! \mathcal{V}  \label{HTR6}
\end{eqnarray}
where \eqref{HTR1} and \eqref{HTR2_1} restricts the vessel pre-dispatch  and trip activation, respectively, subject to purchase decisions  ($a_v^{\mathrm{res}}$). 
\eqref{HTR2} limits the quantity of vessel allocation per load island within $\overline{N}_{s}^{\mathrm{ves}}$. 
Sequential trip execution is imposed in \eqref{HTR3}.
\eqref{HTR4} and \eqref{HTR5} constrain vessel movement to activated trips with assigned decisions. 
Vessel route balance  is constrained in \eqref{HTR6}. 

Moreover, let binary variables $e_{s,v,t}^k$ and $l_{s,v,t}^k$ denote the arrival and  departure  events, while $tr_{s,s',v}$ and $de_{s,v}$ 
represent the transit  and berthing time among islands, respectively. The following constraints define vessels'  temporal consistency with routes. 
\begin{eqnarray}
	&&\textstyle \sum\nolimits_{t\in\mathcal{T}}e_{s,v,t}^k \leq z^k_v, \ \forall s \in\mathcal{S}, \forall v \in \mathcal{V}, \forall k\in\mathcal{K} \label{TE1} \\
	&& \textstyle\sum\nolimits_{s'\in\mathcal{S}} y_{s,s',v}^k( \sum\nolimits_{t\in\mathcal{T}}e_{s,v,t}^k - \sum\nolimits_{t\in\mathcal{T}}l_{s,v,t}^k)=0,\nonumber\\
	&&\qquad\qquad\qquad\qquad \ \forall s \in\mathcal{S}, \forall v \in \mathcal{V}, \forall k\in\mathcal{K} \label{TE2} \\
	&& \textstyle\sum\nolimits_{s'\in\mathcal{S}} y_{s,s',v}^k(\sum\nolimits_{t\!\in\!\mathcal{T}}e_{s,v,t}^kt \!+\! de_{s,v} \!-\! \sum\nolimits_{t\!\in\!\mathcal{T}}l_{s,v,t}^kt)=0, \nonumber\\
	&&\qquad\qquad\qquad\qquad\qquad \forall s \!\in\!\mathcal{S}, \forall v \!\in\! \mathcal{V}, \forall k\!\in\!\mathcal{K} \ \label{TE3} \\
	&& \textstyle	y_{s,s',v}^k\left(\sum\nolimits_{t\in\mathcal{T}}l_{s,v,t}^kt +tr_{s,s',v}- \sum\nolimits_{t\in\mathcal{T}}e_{s',v,t}^kt \right) = 0,\nonumber \\ 
	&& \qquad\qquad\qquad\forall s \in \mathcal{S}, \forall s' \in \mathcal{S}\setminus\mathcal{S}_{d},  \forall k \in \mathcal{K}, \forall v \in \mathcal{V}  \label{TE4}\\
	&&\textstyle y_{s,s',v}^k\left(\sum\nolimits_{t\in\mathcal{T}}l_{s,v,t}^kt +tr_{s,s',v} - \sum\nolimits_{t\in\mathcal{T}}e_{s',v,t}^{k +1}t\right) = 0, \nonumber \\ 
	&&\qquad\qquad\qquad \forall s \in \mathcal{S}, \forall s' \in\mathcal{S}_{d}, \forall k \in \mathcal{K}, \forall v \in \mathcal{V}  \label{TE5}
\end{eqnarray}
where the consistency between arrival decisions and activation  is imposed in \eqref{TE1}. The flow balance is maintained  in \eqref{TE2} among islands. Eqs. \eqref{TE3} defines departure time by linking it to the arrival time, accounting for $de_{s,v}$. Eqs. \eqref{TE4} ensure precedence in vessel movement while \eqref{TE5} extends the precedence condition to consecutive trips, subject to route decisions. Based on their schedule (arrival and departure time), 
the spatial-temporal network for vessel operations, represented by binary time-stamped status $\mu_{sv}^t$ (in-port) and $\beta_v^t$ (at-sea),  is constructed in \eqref{ST1}-\eqref{ST3}. 
\begin{eqnarray}
	&&\textstyle\mu_{s,v}^t\!=\!\sum\nolimits_{t'\!=\!1}^{t}\sum\nolimits_{k\!\in \!\mathcal{K}}(e_{s,v,t'}^k \!-\! l_{s,v,t'}^k)\!,\ \forall s \!\in\!\mathcal{S}\!, \!\forall v \!\in\!\mathcal{V}\!,\!\forall t \!\in\!\mathcal{T} \ \label{ST1}\\
	&&\textstyle\beta_v^t\!=a_v^{\mathrm{ves}}-\!\sum\nolimits_{s\in\mathcal{S}}\!\mu_{s,v}^t, \ \forall v\!\in\!\mathcal{V},\forall t\!\in\!\mathcal{T}\quad \label{ST3}
\end{eqnarray}

\vspace{-5pt}
\subsubsection{Energy Transitions for Maritime Hydrogen Transport Network}
Leveraging the planned routes and schedules of hydrogen vessels, we develop a  maritime hydrogen distribution model. The following constraints govern each vessel’s operations.
\setlength{\arraycolsep}{-0.4em}
\begin{eqnarray}
	&&0\leq m_{s,v,t}^{\mathrm{out}}\leq \mu_{sv}^t \overline{m}_v^{out},\ \forall s\in\mathcal{S}_d, \forall v\in\mathcal{V}, \forall t\in\mathcal{T} \qquad\label{ES1} \\
	&&0\leq m_{s,v,t}^{\mathrm{in}}\leq \mu_{s,v}^t \overline{m}_v^{\mathrm{in}},\quad \forall s\in\mathcal{S}_r,\forall v\in\mathcal{V}, \forall t\in\mathcal{T} \label{ES2} \\
	&& \textstyle\underline{SH}_v\leq\!SH_v^{t+1}\!=\! (\!1 \!- \!\rho_{v})SH_v^{t} \!+  \!\!\Delta t(\!\eta_v\!\sum\nolimits_{s\in\mathcal{S}_r}m_{s,v,t}^{\mathrm{in}}\!-\!\nonumber \\
	&&\qquad\textstyle\!\sum\nolimits_{s\in\mathcal{S}_d}{m_{s,v,t}^{\mathrm{out}}}/{\eta_v}-\beta_v^tm_v^{\mathrm{tr}}) \leq\overline{SH}_v, \ 
	\forall v\in\mathcal{V}, \forall t\in\mathcal{T}  \label{ES3} \\
&& \textstyle \!\!\sum_{t\in\mathcal{T}}(\!\eta_v\!\sum\nolimits_{s\in\mathcal{S}_r}m_{s,v,t}^{\mathrm{in}}\!-\!\textstyle\!\sum\nolimits_{s\in\mathcal{S}_d}{m_{s,v,t}^{\mathrm{out}}}\!/\!{\eta_v}\!-\!\beta_v^tm_v^{\mathrm{tr}})\!\Delta_t \! =\! 0, \forall v \!\in\! \mathcal{V} \label{ES5}
\end{eqnarray}
where the hydrogen  filling ($m_{s,v,t}^{\mathrm{in}}$) and unloading ($m_{s,v,t}^{\mathrm{out}}$) for vessel $v$ is limited in \eqref{ES1} and \eqref{ES2}, receptively, conditional on both in-port status and predefined capacity. 
Storage level $SH_v^{t}$ is imposed in \eqref{ES3} to seek the dynamic balance with $m_{s,v,t}^{\mathrm{in}}$, $m_{s,v,t}^{\mathrm{out}}$, and travel consumption  $m_{s,v,t}^{\mathrm{out}}$, accounting for loss efficiency (i.e., $\rho_{v}$ and ${\eta_v}$) and operational limits [$\underline{SH}_v$,$\overline{SH}_v$]. 
Besides, sustainable operations of hydrogen vessels are ensured in \eqref{ES5}. 
To integrate hydrogen transition within onboard and onshore HSTs, the following constraints model the dynamic flow of hydrogen between vessels and potential nodes across islands. 
\begin{eqnarray}
		&&  \textstyle\sum\nolimits_{v\in\mathcal{V}}m_{s,v,t}^{\mathrm{out}} = \sum\nolimits_{j\in\mathcal{N}_s}m^{\mathrm{din}}_{j,t}, \ \forall s \in \mathcal{S}_d, \forall t \in\mathcal{T} \label{HI1} \\
	&&  \textstyle\sum\nolimits_{v\in\mathcal{V}}m_{s,v,t}^{\mathrm{in}} = \sum\nolimits_{j\in\mathcal{N}_s}m^{\mathrm{dout}}_{j,t}, \ \forall s \in \mathcal{S}_r, \forall t \in\mathcal{T} \label{HI2}
\end{eqnarray}

\section{Two-Stage DDU-DRO Model under Endogenous and Exogenous Uncertainties}
\label{TS_DDUDRO}
A two stage DDU-DRO model is develpoed in this section to ensure both economic efficiency and system resilience under endogenous and exogenous uncertaintiess.
\subsection{1st-Stage Problem: Infrastructure Investment with Proactive Reinforcements}
The first stage problem optimizes the long-term capacity planning of HE-based multi-microgrid system to minimize annual capital expenditure ($\Phi_{\mathrm{capex}}$), which comprises the following four components in \eqref{FO1}: infrastructure investments associated with resource islands ($ \mathfrak{C}_{\mathrm{rs}}$), load islands ($\mathfrak{C}_{\mathrm{ds}} $), hydrogen transport vessels ($ \mathfrak{C}_{\mathrm{hs}}$), and grid reinforcement ($ \mathfrak{C}_{\mathrm{ph}}$) The detailed formulations are  provided below: %
\begin{eqnarray}
	&& \min\nolimits_{\bm{\lambda}\in\Lambda}\quad  \Phi_{\mathrm{capex}} = \mathfrak{C}_{\mathrm{rs}} + \mathfrak{C}_{\mathrm{ds}} + \mathfrak{C}_{\mathrm{hs}}+ \mathfrak{C}_{\mathrm{ph}} \label{FO1}\\ 
	&&\textstyle \mathfrak{C}_{\mathrm{rs}}= \sum\nolimits_{s\in\mathcal{S}_r}\sum\nolimits_{j\in\mathcal{N}_s}\sum\nolimits_{o \in \{wt,pv\}}\!\alpha^{\mathrm{res},o}(c_{\mathrm{res},o,j}^{\mathrm{inv}}\overline{P}_{o,j})\nonumber \\ 
	&&\textstyle\qquad  
	+ \sum\nolimits_{s\in\mathcal{S}_r}(\mathfrak{C}_{\mathrm{elz}}^s+ \mathfrak{C}_{\mathrm{hst}}^s+\mathfrak{C}_{\mathrm{bsb}}^s) \label{FO2} \\
	&&\textstyle  \mathfrak{C}_{\mathrm{ds}}= \sum\nolimits_{s\in\mathcal{S}_d}\sum\nolimits_{j\in\mathcal{N}_s}\alpha^{\mathrm{res}}_{pv}( c_{\mathrm{res},pv}^{\mathrm{inv}}\overline{P}_{pv,j}) \nonumber \\ 
	&&\textstyle\qquad+ \sum\nolimits_{s\in\mathcal{S}_d}( \mathfrak{C}_{\mathrm{hst}}^s+\mathfrak{C}_{\mathrm{bsb}}^s+\mathfrak{C}_{\mathrm{fc}}^s) \label{FO3}\\
	&&\textstyle\mathfrak{C}_{\mathfrak{d}}^s = \sum\nolimits_{j\in\mathcal{N}_s}\alpha^{\mathfrak{d}}( c_{\mathfrak{d}}^{\mathrm{inv}} \overline{P}_{\mathfrak{d},j}),
\  \forall \mathfrak{d} \in\{\mathrm{elz},\mathrm{hst},\mathrm{bsb},\mathrm{fc}\} , \forall s \in\mathcal{S} \label{FO4}\\
	&&\textstyle \mathfrak{C}_{\mathrm{hs}} = \sum\nolimits_{v\in\mathcal{V}}\alpha^{\mathrm{ves}}( c_{\mathrm{pur},v}^{\mathrm{ves}} a^{\mathrm{ves}}_v +\sum\nolimits_{k\in\mathcal{K}}c_{\mathrm{cyc},v,k}^{\mathrm{ves}}z_v^k) \label{FO5} \\
	&&\textstyle \mathfrak{C}_{\mathrm{ph}} = \sum\nolimits_{(i,j)\in\mathcal{B}}\alpha^{\mathrm{gd}} c_{i,j}^{\mathrm{gd}} g_{i,j} + \sum\nolimits_{s\in\mathcal{S}}\sum\nolimits_{j\in\mathcal{N}_{s}}\alpha^{\mathrm{pr}} c^{\mathrm{pr}}_jm^{\mathrm{pr}}_{j}  \label{FO6}
\end{eqnarray}

Note that \eqref{FO2}-\eqref{FO4} incorporate 
capacity-dependent cost ($\bm{c}^{\mathrm{inv}}$),  accounting for expenditures that scale with system capacity. \eqref{FO5} presents the summation of the purchase cost ($\bm{c}^{\mathrm{ves}}_{\mathrm{pur}}$) and trip-active fees ($\bm{c}^{\mathrm{ves}}_{\mathrm{cyc}}$) of hydrogen vessels. In addition, \eqref{FO6} is associated with system enhancement costs including  investment coefficients of line protection ($\bm{c}^{\mathrm{gd}}$) and reserve hydrogen ($\bm{c}^{\mathrm{pr}}$). Also, these costs are evaluated by the capital recovery factor $\alpha^{\mathfrak{d}}=\frac{rr(1+rr)^{Y_{\mathfrak{d}}}}{(1+rr)^{Y_{\mathfrak{d}}}-1}, \mathfrak{d} \in\{\mathrm{res},\mathrm{elz},\mathrm{hst},\mathrm{bsb},\mathrm{fc},\mathrm{ves}, \mathrm{gd},\mathrm{pr}\} $, where $rr$ is the discount rage and  $Y_{\mathfrak{d}}$ is the nominal lifetime.
The first stage decision variables $\bm{\lambda}$ encompass installation or acquisition choices $\bm{a}$, capacity expansion decisions $\bm{\overline{P}}$, vessel routing and schedule configurations ($\bm{x},\bm{z},\bm{y},\bm{e}, \bm{l},\bm{\mu},\bm{\beta}$), protective hardening decisions $\bm{g}$, and storage planning $\bm{m}^{\mathrm{pr}}$, subject to the following constraint space.
\begin{eqnarray}
	&&	\Lambda = \{\ \eqref{HTR1}-\eqref{ST3},\nonumber\\ 
	&&P^{\min}_{o,j}a_{o,j}^{\mathrm{res}}\leq \overline{P}_{o,j}\leq P^{\max}_{o,j}a_{o,j}^{\mathrm{res}}, \   \forall o \in\{wt, pv\}, \forall j\in\mathcal{N}_{s}, \forall s \in \mathcal{S}_r \label{RC1} \\
	&&P^{\min}_{\mathfrak{d},j}a_{j}^{\mathrm{elz}} \leq \overline{P}_{\mathrm{elz},j} \leq P^{\max}_{\mathrm{elz},j}a_{j}^{\mathrm{elz}}, \  \forall j\in\mathcal{N}_{s}, \forall s \in\mathcal{S}_r \label{RC2} \\
	&&P^{\min}_{pv,j}a_{pv,j}^{\mathrm{res}}\leq \overline{P}_{pv,j}\leq P^{\max}_{pv,j}a_{pv,j}^{\mathrm{res}}, \  \forall j\in\mathcal{N}_{s}, \forall s \in \mathcal{S}_d \label{LC1} \\
	&&P^{\min}_{\mathfrak{d},j}a_{j}^{\mathfrak{d}} \leq \overline{P}_{\mathfrak{d},j} \leq P^{\max}_{\mathfrak{d},j}a_{j}^{\mathfrak{d}}, \ \mathfrak{d} \in \{\mathrm{hst, bsb}\}, \forall j\in\mathcal{N}_{s}, \forall s\in\mathcal{S} \label{ESS} \\
	&&\textstyle\sum\nolimits_{(i,j)\in \mathcal{B}}g_{i,j} \leq N_{\mathrm{line}}^{\max}  \label{GR1}\\
	&&\textstyle \sum\nolimits_{s\in\mathcal{S}}\sum\nolimits_{j\in\mathcal{N}_s}m^{\mathrm{pre}}_j \leq N_{\mathrm{str}}^{\max} \} \label{GR2}
\end{eqnarray}
where \eqref{RC1} restricts the installed capacity $\overline{P}_{o,j}$ of WT and PV units in resource islands with specified bounds, subject to the installation decision $a_{o,j}^{\mathrm{res}}$. Also,  \eqref{RC2} imposes the capacity limits on the  electrolyzers conditional on their deployment decision $a_{j}^{\mathrm{elz}}$.  As for load islands, \eqref{LC1}  constrains capacity limits on PV systems. Notice that the energy storage capacities are subject to \eqref{ESS}, governing the HSTs and BSBs in all islands. 
Besides,  \eqref{GR1} limits the maximum protective actions and \eqref{GR2} enforces the upper bound of the pre-installed storage. 
\subsection{Endogenous and Exogenous Uncertainties}
\subsubsection{Wind-Induced  Uncertainty with Contingency Events}
\begin{figure}[!t]
	\vspace{-10pt}
	\centering
	\includegraphics[width=4in]{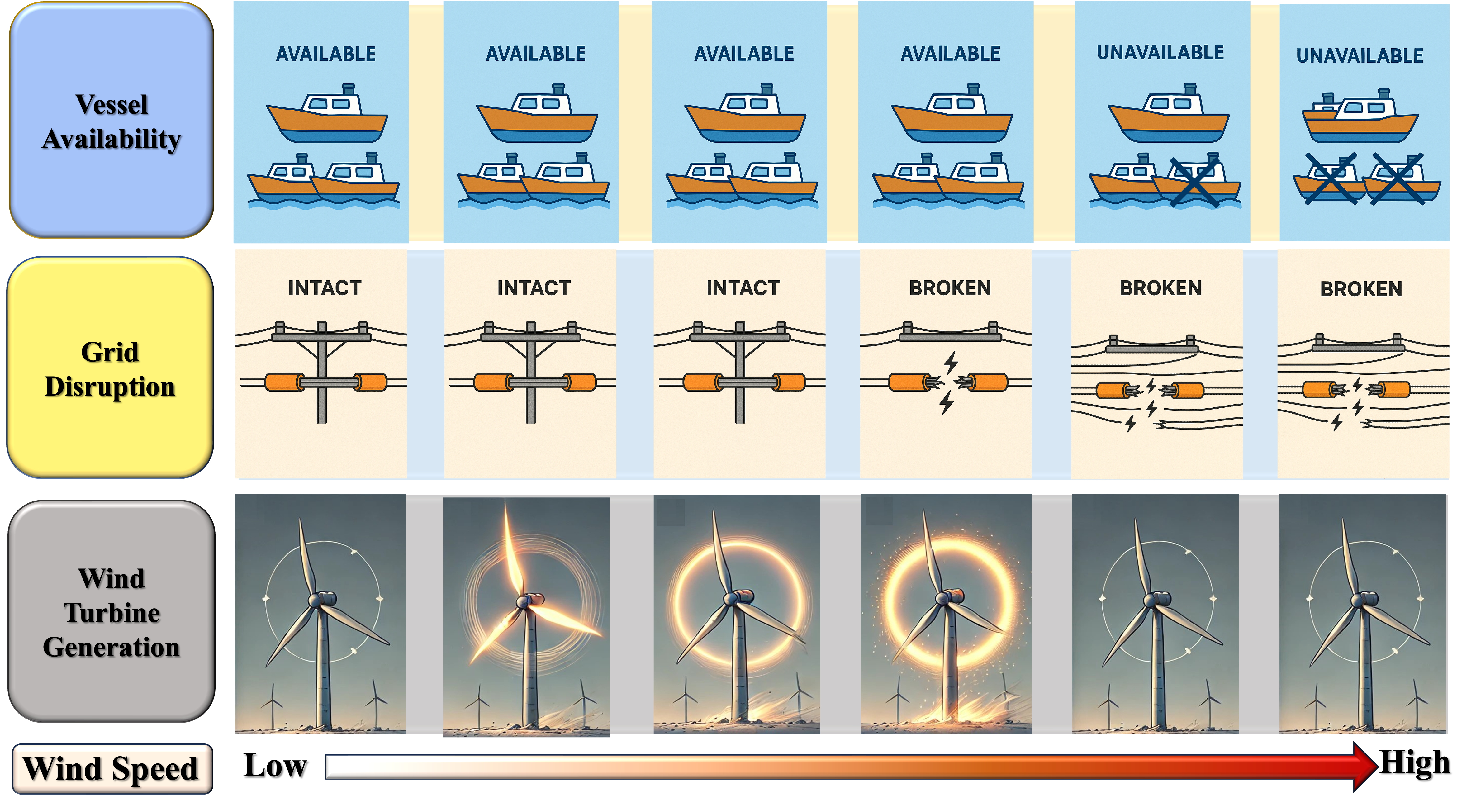}
	\caption{Impact of varying wind speed levels  on vessel availability, grid disruption, and wind power generation} 
	\label{fig:uncertain} 
\end{figure}

Wind is the most influential factor in island multi-microgrid systems, serving as both a critical RES and a major constraint. Specifically, at low wind speeds, insufficient airflow prevents turbines from generating power while moderate winds allow stable power generation. However, excessively strong winds can trigger turbine shutdowns for safety reasons, which also disrupt grid stability and hinder the operation of hydrogen  vessels. The impact of wind variations is shown in Fig. \ref{fig:uncertain}. 
To systematically capture its impacts, wind speed is categorized into discrete levels $u\in\mathcal{U}$, each representing a bounded range of velocities, whose realizations (denoted as random binary  variable $\upsilon_u$) are constrained by \eqref{U1}. 
\begin{eqnarray}
		&&   \textstyle\sum\nolimits_{u \in \mathcal{U}}\upsilon_u = 1 \label{U1}
\end{eqnarray} 

Note that wind energy utilization factor ($\tilde{\delta}_{wt, j}^t$) follows a nonlinear relationship with wind speed. 
To model this, leveraging historical data and turbine power curve,  experts can  estimate $\tilde{\delta}_{wt, j}^t$ by \eqref{F1}, using the nominal wind power output $\overline{\delta}_{wt, j}^{u, t}$, along with the maximum upward ($\hat{\delta}^{u,t}_{wt, j}$) and downward ($\check{\delta}^{u,t}_{wt, j}$) deviations. 
\begin{eqnarray}
	&& \textstyle\tilde{\delta}_{wt,j}^t \! = \!  \sum\limits_{u \in \mathcal{U}}\upsilon_u \overline{\delta}_{wt, j}^{u, t}\! +\!  \tau_{wt, j}^{t, +}\sum\limits_{u \in \mathcal{U}}\upsilon_u \hat{\delta}^{u,t}_{wt, j} \! -\!  \tau_{wt, j}^{t, -}\sum\limits_{u \in \mathcal{U}}\upsilon_u\check{\delta}^{u,t}_{wt, j} \nonumber\\
	&&\qquad\qquad\qquad\qquad\qquad \forall j \in\mathcal{N}_s, \forall s \in\mathcal{S}_r, \forall t \in \mathcal{T} \label{F1}
\end{eqnarray}
where $\tau_{wt, j}^{t, +}$ and $\tau_{wt, j}^{t, -}$, within [0, 1], represent the upward and downward deviation factors, respectively.

Then, wind speed can pose significant risks to both power grids and maritime hydrogen transport networks, especially as  $N-k$ contingencies become more severe under extreme conditions.
\begin{eqnarray}
	&&   \textstyle\sum\nolimits_{(i,j)\in \mathcal{B}}(1-\tilde{\omega}_{ij})\leq \sum\nolimits_{u \in \mathcal{U}}\upsilon_uD_u^{\mathrm{line}} \\
	&&   \textstyle\sum\nolimits_{v\in\mathcal{V}}(1-\tilde{\phi}_{v})\leq \sum\nolimits_{u \in \mathcal{U}}\upsilon_uD_u^{\mathrm{ves}}  
\end{eqnarray}
where $\tilde{\omega}_{ij}$ and  $\tilde{\phi}_{v}$ are binary variables indicating whether transmission line $(i,j)$ and vessel $v$ remain intact (1) or have failed (0). The uncertainty budgets $D_u^{line}$ and $D_u^{vss}$ constrain the maximum number of disrupted lines and vessels for each wind level $u$. 

\subsubsection{Other Exogenous Uncertainty}Solar power generation and electricity nodal active/reactive demand  in load islands are treated as independent uncertainties, each with its own constraints and ranges, as shown in the following \eqref{IU1}-\eqref{IU3}.
\begin{eqnarray}
	&&  \underline{\delta}_{pv, j, t} \leq\tilde{\delta}_{pv,j}^t \leq  \overline{\delta}_{pv, j, t},  \  \forall j\in\mathcal{N}_{s}, \forall s \in\mathcal{S},\forall t\in\mathcal{T}\label{IU1} \\
	&&\underline{p}_{\mathrm{l},j, t} \leq \tilde{p}_{\mathrm{l},j}^t \leq  \overline{p}_{\mathrm{l},j,t},\  \forall j\in \mathcal{N}^a_s, \forall s \in\mathcal{S}_d, \forall t\in\mathcal{T}  \label{IU2} \\
	&& \underline{q}_{\mathrm{l},j, t} \leq \tilde{q}_{\mathrm{l},j}^t \leq \overline{q}_{\mathrm{l},j,t}, \  \forall j\in \mathcal{N}^a_s,\forall s\in\mathcal{S}_d, \forall t\in\mathcal{T}   \label{IU3} 
\end{eqnarray}

\subsubsection{Endogenous Uncertainty}We note that constraints \eqref{IU4}–\eqref{IU5} represent decision-dependent uncertainties, where the occurrence of the uncertain active/reactive power demand depend on indicator variables  $\iota_{j}^t =\mathbbm{1}(\sum_{o\in \{wt,pv\}}a_{o,j}^{\mathrm{res}}\geq1)$. The uncertainty realizations are only triggered if at least one RES is allocated, reflecting the auxiliary demand arising from associated renewable energy generation and hydrogen production. 
\begin{eqnarray}
		&&\iota_{j}^t \underline{p}_{\mathrm{l},j, t} \leq \tilde{p}_{\mathrm{l},j}^t \leq  \iota_{j}^t \overline{p}_{\mathrm{l},j,t},\  \forall j\in \mathcal{N}_s, \forall s\in\mathcal{S}_r, \forall t\in\mathcal{T}  \label{IU4} \\
	&& \iota_{j}^t \underline{q}_{\mathrm{l},j, t} \leq \tilde{q}_{\mathrm{l},j}^t \leq\iota_{j}^t  \overline{q}_{\mathrm{l},j,t}, \  \forall j\in \mathcal{N}_s,\forall s\in\mathcal{S}_r, \forall t\in\mathcal{T}   \label{IU5} 
\end{eqnarray}

Thus, the set of all possible uncertainty realizations can be presented by DDU-based sample space $\Omega(\bm\lambda) =\{\bm{\xi} \subseteq\mathbb{R}^{n_1}\times\{0, 1\}^{n_2}|\eqref{U1}-\eqref{IU5}\}$, where $\bm{\xi} =\{\bm{\upsilon},\bm{\tau}^{+},\bm{\tau}^{-},\bm{\tilde{\omega}}, \bm{\tilde{\phi}}, \bm{\tilde{\delta}}, \bm{p}, \bm{q}\}$. 
\subsubsection{DDU-based Ambiguity Set}Rather than considering a single and precise distribution for these uncertainties, we define the DDU-based ambiguity set $\mathcal{P}(\bm\lambda)=\{\mathbb{P}\in\mathcal{M}(\Omega(\bm{\lambda}), \mathcal{F})|\eqref{UM0}-\eqref{UM2}\}$ where  $\mathcal{F}$ is $\sigma$-algebra containing all singletons and \eqref{UM0}-\eqref{UM2} are first-order moment constraints. 
\begin{eqnarray}
&& \mathbb{E}_{\mathbb{P}}[1]=1 \label{UM0}\\
&&	\underline{\upsilon}_u \leq \mathbb{E}_{\mathbb{P}}[\upsilon_u] \leq \overline{\upsilon}_u, \  \forall u \in\mathcal{U} \label{UM1} \\
&&	\underline{\phi}_v \leq \mathbb{E}_{\mathbb{P}}[\tilde{\phi}_{v}] \leq \overline{\phi}_v, \  \forall v \in\mathcal{V} \label{UM3}\\
&&	\underline{\mathfrak{y}}_{j,t}^\mathrm{e} \leq \mathbb{E}_{\mathbb{P}}[\mathfrak{y}_{j}^t] \leq \overline{\mathfrak{y}}_{j, t}^{\mathrm{e}}, \mathfrak{y}\in\{\tau_{wt}^{+},\tau_{wt}^{-},\tilde{\delta}_{pv},\tilde{p}_{\mathrm{l}}, \tilde{q}_\mathrm{l}\} \nonumber\\
&&\qquad\qquad\qquad\qquad\qquad\forall j \in\mathcal{N}_s, \forall s \in\mathcal{S}, \forall t \in\mathcal{T}  \label{UM4}\\
&&	\underline{\omega}_{i,j} \leq \mathbb{E}_{\mathbb{P}}[\tilde{\omega}_{i,j}] \leq \overline{\omega}_{i,j}-\varepsilon_{ij}g_{i,j}, \forall (i, j) \in\mathcal{B}\label{UM2} 
\end{eqnarray}

Notably, \eqref{UM2} introduces a distinct form of decision-dependent moment-based constraints by adjusting the upper bound of $ \mathbb{E}_{\mathbb{P}}[\tilde{\omega}_{i,j}]$. Here, $\varepsilon_{ij}$ is a reduction factor reflecting the extent to which the grid hardening decision $g_{i,j}$ narrows the allowable expectation range.

\subsection{2nd--Stage Problem: Energy Dispatch and Cross-Island Delivery}
As shown in \eqref{SeT1}, the second stage aims to determine the worst-case probability distribution within the DDU consideration. Given the first-stage investment decisions $\bm{\lambda}$ and the realization of uncertainty $\bm{\xi}$, the recourse function $Q(\bm{\lambda},\bm{\xi})$ seeks to minimize operational and maintenance (O\&M) costs. These costs, as formulated  in  \eqref{SeT2},  associate with renewable hydrogen production in resource islands ($\mathfrak{E}^{\mathrm{o\&m}}_{\mathrm{rs}, t}$), electricity system operation in load islands ($\mathfrak{E}^{\mathrm{o\&m}}_{\mathrm{ds}, t}$),  and maritime hydrogen transportation ($\mathfrak{E}^{\mathrm{o\&m}}_{\mathrm{hs},t}$), all of which are annualized through a scaling factor $\hbar$. The detailed formulations are as follows: 
\begin{eqnarray}
	&&\sup\nolimits_{\mathbb{P}\in\mathcal{P}(\bm{\lambda})} \quad \mathbb{E}_{\mathbb{P}}[Q(\bm{\lambda},\bm{\xi})] \label{SeT1} \\ 
	&& \textstyle Q(\bm{\lambda},\bm{\xi}) = \min\nolimits_{\bm{\zeta}\in\mathcal{Z}(\bm{\lambda},\bm{\xi})}  \hbar\sum\nolimits_{t\in\mathcal{T}} (\mathfrak{E}^{\mathrm{o\&m}}_{\mathrm{rs}, t} + \mathfrak{E}^{\mathrm{o\&m}}_{\mathrm{ds}, t} +\mathfrak{E}^{\mathrm{o\&m}}_{\mathrm{hs},t}) \label{SeT2} \\
	&&\textstyle\mathfrak{E}^{\mathrm{o\&m}}_{\mathrm{rs}, t} \! =\! \sum\nolimits_{s\in\mathcal{S}_r} \sum\nolimits_{j\in\mathcal{N}_{s}}(\sum\nolimits_{o\in \{wt,pv\}} O_o^{\mathrm{res}} p^{t}_{o,j}+\nonumber\\
	&&\textstyle\qquad\quad O^{\mathrm{hst}}m^{\mathrm{dout}}_{j,t} + O^{\mathrm{bsb}}(p_{\mathrm{ch},j}^{t}+p_{\mathrm{dis},j}^{t} ), \ \forall t\in\mathcal{T}  \label{SeT3}\\ 
	&&\mathfrak{E}^{\mathrm{o\&m}}_{\mathrm{ds}, t} =  \textstyle\sum\nolimits_{s\in\mathcal{S}_d}\sum\nolimits_{j\in\mathcal{N}_{s}}( O^{\mathrm{res}}_{pv} p^{t}_{pv,j}+O^{\mathrm{hst}}m_{j,t}^{\mathrm{din}}+ \nonumber\\
	&&\qquad\quad O^{\mathrm{bsb}}(p_{\mathrm{ch},j}^{t}+p_{\mathrm{dis},j}^{t} ), \ \forall t\in \mathcal{T} \label{SeT4}\\ 
	&&\textstyle\mathfrak{E}^{\mathrm{o\&m}}_{\mathrm{hs}, t}  =\textstyle\sum\nolimits_{s\in\mathcal{S}}\sum\nolimits_{v\in\mathcal{V}}O^{\mathrm{ves}}\left(m^{\mathrm{in}}_{s,v,t}+m^{\mathrm{out}}_{s,v,t} \right), \ \forall t \in\mathcal{T} \label{SeT5} 
\end{eqnarray}
The above cost terms are calculated by the operational cost coefficients $\bm{O}$ associated with different system components. The second-stage decision variables $\bm\zeta$ , include the hydrogen production and storage at resource islands, the amount of  hydrogen transported energy between islands, and the hydrogen-electrical system's operations in load islands. These variables must satisfy the following constraints $\mathcal{Z}(\bm{\lambda},\bm{\xi})$:
\begin{eqnarray}
&&\mathcal{Z}(\bm{\lambda},\bm{\xi}) = \{ \eqref{RE1} -\eqref{PDN8}, \eqref{ES1} -\eqref{HI2}, \nonumber \\ 
	&&\textstyle| \sum\nolimits_{t\in\mathcal{T}}(\eta_{\mathrm{hst+}}m_{\mathrm{elz},j}^t-{m_{j,t}^{\mathrm{dout}}}/{\eta_{\mathrm{hst-}}} )| \leq m^{\mathrm{pr}}_j, \  \forall j\in\mathcal{N}_{s}, \forall s\in\mathcal{S}_r  \label{SS1}\\ 
	&&\textstyle| \sum\nolimits_{t\in\mathcal{T}}(\eta_{\mathrm{hst+}} m_{j,t}^{\mathrm{din}}-{m_{\mathrm{fc},j}^t}/{\eta_{\mathrm{hst-}}} )|\leq m^{\mathrm{pr}}_j, \  \forall j\in\mathcal{N}_{s}, \forall s\in\mathcal{S}_d\label{SS2}\\
	&&\tilde{\omega}_{i,j}(U_i^t\!-\!U_j^t\! -\! ({r_{i,j}fp_{i,j}^t \! +\! x_{i,j}fq_{i,j}^t})/{U_0})\! =\! 0, \forall (i,j)\in \mathcal{E}, \forall t\in\mathcal{T}  \label{SS3} \\
	&&-\tilde{\omega}_{i,j}\overline{FP}_{ij}\leq fp_{ij}^t\leq \tilde{\omega}_{i,j}\overline{FP}_{ij}, \ \forall (i,j)\in \mathcal{E}, \forall t\in\mathcal{T} \label{SS4}\\
	&& \textstyle 0 \leq p_{\mathrm{l},j}^t =  \tilde{p}_{\mathrm{l},j}^t - p_{\mathrm{ls},j,t}, \ \forall j\in \mathcal{N}, \forall t\in\mathcal{T}  \label{SSB1} \\
	&& \textstyle \sum_{j\in\mathcal{N}}p_{\mathrm{ls}, j,t} \leq \sum_{u \in \mathcal{U}}\upsilon_u  \overline{P}_{\mathrm{ls},j,t}^{u}, \  \forall t\in\mathcal{T} \label{SSB2} \\
	&& \textstyle 0\leq q_{\mathrm{l},j}^t =  \tilde{q}_{\mathrm{l},j}^t - q_{\mathrm{ls},j,t}, \ \forall j\in \mathcal{N}, \forall s \in\mathcal{S},\forall t\in\mathcal{T}  \label{SSB3} \\
	&& \textstyle \sum_{j\in\mathcal{N}}q_{\mathrm{ls}, j,t} \leq \sum_{u \in \mathcal{U}}\upsilon_u  \overline{Q}_{\mathrm{ls},j,t}^{u}, \  \forall t\in\mathcal{T} \label{SSB4} \\
	&&0\leq m_{svt}^{\mathrm{out}}\leq(1- \tilde{\phi}_{v})\mu_{sv}^t \overline{m}_v^{\mathrm{out}},\ \forall s\in\mathcal{S}_d, \forall m\in\mathcal{V}, \forall t\in\mathcal{T} \label{SS5} \\
	&&0\leq m_{svt}^{\mathrm{in}}\leq( 1- \tilde{\phi}_{v})\mu_{sv}^t \overline{m}_v^{\mathrm{in}},\ \forall s\in\mathcal{S}_r,\forall v\in\mathcal{V}, \forall t\in\mathcal{T} \} \label{SS6} 
\end{eqnarray}
where constraints \eqref{SS1} and \eqref{SS2} govern the adjustable balance of HSTs at resource and load islands over a  period, which can mitigate the impact  uncertainties by allowing pre-assigned acceptable hydrogen buffering storage  ($\bm{m}^{\mathrm{pr}}$). 
\eqref{SS3} and \eqref{SS4} illustrate the relationship between the power flow and line status  ($\tilde{\omega}_{i,j}$).
Besides, $\mathcal{N}$ denotes the set of nodes involved in the system, where $\mathcal{N} := \left( \bigcup_{s \in \mathcal{S}_d} \mathcal{N}_s^a \right) \cup \left( \bigcup_{s \in \mathcal{S}_r} \mathcal{N}_s \right)$. Let $p_{\mathrm{ls},j,t}$ and $q_{\mathrm{ls},j,t}$ denote the load shedding at each node in \eqref{SSB1} and \eqref{SSB3}, with their total shedding constrained by the maximum shedding capacity for active and reactive power under wind level $u$ in \eqref{SSB2} and \eqref{SSB4}, respectively. Also, the maritime hydrogen transport network is subject to constraints \eqref{SS5} and \eqref{SS6} that limit hydrogen unloading and refilling rates, particularly in the presence of vessel failures.

\subsection{Compact Formulation}
To enhance conciseness while preserving generality, we reformulate the above problem as the following $\mathbf{DDU-DRO}$ model in a compact matrix representation:
\begin{eqnarray}
	&&\mathbf{DDU-DRO}: \ \min_{\bm{\lambda}} \ \bm{c}^\mathrm{T}\bm{\lambda}+\sup_{\mathbb{P}\in\mathcal{P}(\bm \lambda)} \int_{\bm{\xi} \in\Omega(\bm\lambda)}Q(\bm{\lambda},\bm{\xi})\mathbb{P}(d\bm{\xi}) \label{CP1} \\ 
	&&\mathrm{s.t.} \ \  \Lambda=\left\lbrace \bm{\lambda}\in \{0,1\}^{n_1}\times \mathbb{R}_+^{n_2} : \bm{A}\bm{\lambda}\leq \bm{b} \right\rbrace   \label{CP2}\\
	&& \qquad\mathcal{P}(\bm\lambda) \!= \!\{\mathbb{P} \!\in\!\mathcal{M}(\Omega(\bm{\lambda}), \mathcal{F})\!:\!\mathbb{E}_{\mathbb{P}}[1]\!=\!1,\mathbb{E}_{\mathbb{P}}[\bm{\psi}(\bm{\xi})]\!\leq\!\bm{\gamma}_1 - \bm{\gamma}_2\bm{\lambda}\} \label{CP3} \\ 
	&&\qquad\Omega(\bm\lambda) = \{\bm{\xi}\in \{0,1\}^{n_3}\times \mathbb{R}_+^{n_4}: \bm{H}\bm{\xi} \leq\bm{k}_1 , \bm{\xi} \leq \bm{k}_2\bm{\lambda} \} \label{CP4}\\
	&&\qquad Q(\bm{\lambda},\bm{\xi}) = \min_{\bm{\zeta}}  \bm{g}^{\mathrm{T}}\bm{\zeta}  \label{CP5}\\ 
	&&\qquad\mathcal{Z}(\bm{\lambda},\bm{\xi})=\{ \bm{\zeta}\in \mathbb{R}^{n_5} :  \bm{F} \bm{\zeta} \geq \bm{f}-\bm{L} \bm{x} - \bm{G}\bm{\xi}\} \label{CP6} 
\end{eqnarray} 

The objective function in \eqref{CP1} consists of two components: the first-stage investment cost and the worst-case expectation of the second-stage O\&M cost, corresponding to \eqref{FO1}-\eqref{FO6} and \eqref{SeT1}-\eqref{SeT5}, respectively. The first-stage decision variables, denoted by vector $\bm{\lambda}$, are constrained within the feasible set $\Lambda$, as defined in \eqref{CP2}. The probability distribution is characterized by the ambiguity set $\mathcal{P}(\bm \lambda)$, which is constructed based on moment inequalities, as specified in \eqref{CP3}. The sample space of random parameters is denoted by $\Omega(\bm \lambda)$, with feasibility constraints given in \eqref{CP4}. The second-stage recourse problem, represented by \eqref{CP5}, determines the optimal operational adjustments to mitigate the impact of uncertainty. $\mathcal{Z}(\bm{\lambda},\bm{\xi})$ in \eqref{CP6} is the feasible set of the corresponding decision variables $\bm{\zeta}$. These constraints encapsulate the recourse actions under the realization of $\bm{\lambda}$ and $\bm{\xi}$, incorporating both operational and system constraints.
\section{A Primal-based Decomposition Algorithm for DRO Model with Strong Cuting Planes}
\label{PDDRO}
The proposed $\mathbf{DDU-DRO}$ model presents significant computational challenges due to its tri-level structure, the ambiguity set that contains uncountable many probability distributions,  and the complexity introduced by the DDU-based sample space. Unlike previous studies, where recourse problem is assumed to be always feasible for any first-stage decision, our model encounters infeasibility issue in the recourse problem under certain realizations of random factors, which dramatically increases the model's complexity.  To efficiently address this issue,  we develop strong cutting planes \cite{lu2024two}, which can be directly embedded into a classical C\&CG framework\cite{zeng2013solving}. From the primal perspective, we further design a column generation (CG) procedure to solve the inner bilevel worst-case expectation problem. %
\subsection{Neutralization Reformulation of DDU-based  Sample Space}
The DDU-based sample space $\Omega(\bm\lambda)$ induces a non-convex and combinatorially complex computational challenges, which substantially complicates the tractability of the proposed $\mathbf{DDU-DRO}$ formulation. To address this issue, we adopt a neutralization reformulation approach that transforms sample space from the original DDU structure $\Omega(\bm{\lambda})$ into a tractable decision independent uncertainty (DIU) representation $\Omega$\cite{zeng2022two}. 
\begin{prop}
	\label{DDU}
	Formulation $\mathbf{DDU-DRO}$ is equivalent to the following  $\mathbf{DDAS-DRO}$ with DIU-based sample space:
	\begin{eqnarray}
		&&\mathbf{DDAS-DRO}: \ \min_{\bm{\lambda}\in\Lambda} \ \bm{c}^\mathrm{T}\bm{\lambda}+\sup_{\mathbb{P}\in\mathcal{P}(\bm\lambda)} \int_{\bm{\xi} \in\Omega}Q(\bm{\lambda},\bm{\xi})\mathbb{P}(d\bm{\xi})
	\end{eqnarray}
	where
	\begin{eqnarray}
		&&\Omega = \{\bm{\xi}\in \{0,1\}^{n_3}\times \mathbb{R}_+^{n_4}: \bm{H}\bm{\xi} \leq\bm{k}_1\} \label{CP4_1}\\
		&& Q(\bm{\lambda},\bm{\xi}) = \min_{\bm{\zeta}}  \bm{g}^{\mathrm{T}}\bm{\zeta}  \label{CP5_1}\\ 
		&&\mathcal{Z}(\bm{\lambda},\bm{\xi})=\big\{ \bm{\zeta}\in \mathbb{R}^{n_5} :  \bm{F} \bm{\zeta} \geq \bm{f}-\bm{L} \bm{x} - \bm{G}\bm{\xi}\ - \bm{k}_3\big(\bm{\xi}\circ(\bm{{k}_2\lambda})\big)\big\} \  \label{CP6_1} 
	\end{eqnarray}
\end{prop}

Note that $\bm{\xi} \circ (\bm{k}_2\bm{\lambda})$ is a Hadamard-type projection operator to neutralize infeasible scenarios (i.e., those in $\Omega \setminus \Omega(\bm\lambda)$), where the operator $\circ$ denotes element-wise multiplication. This operation effectively maps any infeasible scenario $\bm{\xi} \notin \Omega(\bm\lambda)$ into a feasible one by zeroing out components forbidden by the current first-stage decision $\bm{\lambda}$. Based on Proposition \ref{DDU}, we preserve the decision-dependent nature of the original model while enabling tractable reformulations $\mathbf{DDAS-DRO}$ under the DIU-based sample space.
\subsection{Discrete Reformulation of the Worst-Case Expected Value } 
The $\mathbf{WCEV}$ problem, i.e, $\sup_{\mathbb{P}\in\mathcal{P}(\bm \lambda)}\mathbb{E}_{\mathbb{P}}[{Q}(\bm{\lambda},\bm{\xi})]$, involves optimizing under the worst-case probability distribution, which is inherently intractable due to its integral formulation over an ambiguity set with uncountable many probability distributions. 
To address this complexity, we reformulate the second stage problem into an equivalent discrete representation by leveraging Proposition \ref{p1}.
\begin{prop}
	\label{p1}
	For a given $\bm \lambda$, let $\{(\bm{\xi}_j, p_{\bm{\xi}_j})\}_{j=1}^n$ denote a discrete probability distribution in the ambiguity set $\mathcal{P}(\bm \lambda)$, where each uncertain scenario $\bm{\xi}_j$ is associated with its probability $p_{\bm{\xi}_j}$. Given that the functions $Q(\bm{\lambda}, \bm{\xi})$ and $\bm{\psi}(\bm{\xi})$ are continuous over $\bm \xi \in \Omega$, the following equivalence holds:
	\begin{eqnarray}
		&&\sup_{\mathbb{P}\in\mathcal{P}(\bm \lambda)} \int_{\bm{\xi} \in\Omega}Q(\bm{\lambda},\bm{\xi})\mathbb{P}(d\bm{\xi}) = \max_{\mathbb{P}\in\mathcal{P}(\bm \lambda)}  \int_{\bm{\xi} \in\Omega}Q(\bm{\lambda},\bm{\xi})\mathbb{P}(d\bm{\xi}) \nonumber \\
		&&= \lim_{n\rightarrow+\infty}\max_{\{(\bm{\bm{\xi}}_j,p_{{\bm{\xi}}_j})\}_{j=1}^n\in\mathcal{P}(\bm \lambda)} \sum_{j=1}^nQ(\bm{\lambda},\bm{\xi}_j)p_{{\bm{\xi}_j}}  \\ 
		&& \geq \max_{p_{{\bm{\xi}}^{\mathrm{f}}}} \big\{\sum_{j=1}^nQ(\bm{\lambda},\bm{\xi}_j^{\mathrm{f}})p_{{\bm{\xi}_j^{\mathrm{f}}}}: \{(\bm{\bm{\xi}}^{\mathrm{f}}_j,p_{{\bm{\xi}}^{\mathrm{f}}_j})\}_{j=1}^n \in \mathcal P(\bm \lambda)\big\} \nonumber
	\end{eqnarray}
	where $\{{\bm{\xi}}_j^{\mathrm{f}}\}_{j=1}^n$ is a set of fixed scenarios.
\end{prop} 

Notice that Proposition \ref{p1} not only converts the original integral-based formulation into an infinite-dimensional summation problem, but also provides a lower bound for the optimal value of $\mathbf{WCEV}$, which can be tightened by including additional nontrivial scenarios.

\subsection{Strong Cutting Planes: Addressing Intractable Feasibility Barriers}
Most existing studies implicitly or explicitly adopt the so-called complete recourse assumption, i.e., it is presumed that   $Q(\bm\lambda, \bm{\xi})$ is always feasible for any  $\bm\lambda$ and $\bm\xi$. Based on this assumption and Proposition \ref{p1}, the main master problem (denoted as $\mathbf{MMP}$), serving  as a relaxation of the primal DRO model,  is  formulated as follows: 
\begin{eqnarray}
	&&\textstyle \mathbf{MMP}: \ \Gamma = \min_{\bm{\lambda}\in\Lambda} \ \bm{c}^\mathrm{T}\bm{\lambda}+\eta \label{MMP1} \\ 
	&&\textstyle \mathrm{s.t.}\ \  \eta \geq \max\{\sum_{\bm\xi\in \Xi}(\bm{g}^\mathrm{T}\bm{\zeta}_{\bm\xi})p_{\bm{\xi}}: (p_{\bm\xi})_{\bm\xi\in\Xi}\in\mathcal{P}(\bm \lambda) \} \label{MMP2-1}\\ 
	&&\textstyle  \qquad\bm\zeta_{\bm\xi}\in {\mathcal{Z}}(\bm\lambda,\bm{\xi}), \ \forall \bm{\xi}\in\Xi _\mathrm{o}\label{MMP4}
\end{eqnarray} 
where $\Xi_\mathrm{o}$ is the fixed set of scenarios, which could be  provided by  $\mathbf{WCEV}$. However, such an assumption 
is highly restrictive in energy system applications, especially under deep uncertainty or stressed system conditions such as post-contingency or high renewable penetration.
To address this challenge, we introduce a feasibility-checking procedure via a designed feasibility problem $\mathbf{WCEV}(F)$: $\sup_{\mathbb{P}\in\mathcal{P}(\bm \lambda)}\mathbb{E}_{\mathbb{P}}[\tilde{Q}(\bm{\lambda},\bm{\xi})]$ as follows: 
\begin{eqnarray}
	&&\textstyle	\tilde{Q}(\bm{\lambda}, \bm{\xi}) = \min_{\bm{\zeta}\in\tilde{\mathcal{Z}}} ( \max\{0,  \sum_{t\in\mathcal{T}}(( \sum_{j\in\mathcal{N}}p_{\mathrm{ls}, j, t} - \nonumber \\
	&&\textstyle \sum_{u \in \mathcal{U}}\upsilon_{ij}\overline{P}_{\mathrm{ls}, j, t}) \!+\!(\sum_{j\in\mathcal{N}}q_{\mathrm{ls}, j, t} \!-\!\sum_{u \in \mathcal{U}}\upsilon_{ij}\overline{Q}_{\mathrm{ls}, j, t} )) \})
\end{eqnarray}
whose compact formulation is as $\tilde{Q}(\bm{\lambda}, \bm{\xi}) =  \min_{\bm{\zeta}\in\tilde{\mathcal{Z}}(\bm{\lambda},\bm{\xi})}  \bm{\kappa}^{\mathrm{T}}\bm{\zeta}$ with  $\tilde{\mathcal{Z}}(\bm{\lambda},\bm{\xi}) = \{\mathcal{Z}\setminus\{\eqref{SSB2},\eqref{SSB4}\}\}$. According to the form of $\mathbf{WCEV}(F)$, the following result holds directly. %
\begin{prop}
	\label{FR}
	A first-stage decision $\bm{\lambda}$ is almost surely feasible, i.e., the infeasible scenarios occur  with zero probability, for the recourse problem $Q(\bm{\lambda}, \bm{\xi})$ if and only if
	$\max\nolimits_{\mathbb{P}\!\in\!\mathcal{P}(\bm \lambda)}\mathbb{E}_{\mathbb{P}}[\tilde{Q}(\bm{\lambda},\bm{\xi})]\! =\! 0$.
	Accordingly, the $\mathbf{DDAS-DRO}$ problem is infeasible if no such $\bm{\lambda}$ exists.
\end{prop}
Proposition~\ref{FR} provides a necessary and sufficient condition for verifying the feasibility of the first-stage decisions. This condition, though implicit, should be explicitly incorporated into the  $\mathbf{MMP}$ without sacrificing optimality so that the feasible region of $\bm\lambda$ only consists of almost sure feasible solutions. Hence, we introduce strong cutting planes, as shown in the revised $\mathbf{MMP}$ formulation below. 
\begin{eqnarray}
	&&\textstyle \mathbf{MMP}: \ \Gamma = \min_{\bm{\lambda}\in\Lambda} \ \bm{c}^\mathrm{T}\bm{\lambda}+\eta \label{MMP1} \\ 
	&&\textstyle \mathrm{s.t.}\ \  \eta \geq \max\{\sum_{\bm\xi\in \Xi}(\bm{g}^\mathrm{T}\bm{\zeta}_{\bm\xi})p_{\bm{\xi}}: (p_{\bm\xi})_{\bm\xi\in\Xi}\in\mathcal{P}(\bm \lambda) \} \label{MMP2}\\ 
	&& \textstyle\qquad0 \geq \max\{\sum_{\bm\xi\in \Xi}(\bm{\kappa}^\mathrm{T}\bm{\zeta}_{\bm\xi})p'_{\bm{\xi}}:(p'_{\bm\xi})_{\bm\xi\in\Xi}\in\mathcal{P}(\bm \lambda)\} \label{MMP3}\\
	&&\textstyle  \qquad\bm\zeta_{\bm\xi}\in\tilde{\mathcal{Z}}(\bm\lambda,\bm{\xi}), \ \forall \bm{\xi}\in\Xi \label{MMP4}
\end{eqnarray} 
where $\Xi$ is the fixed set of scenarios, which could be provided by the $\mathbf{WCEV}$ and $\mathbf{WCEV}(F)$. \eqref{MMP2} and \eqref{MMP3} represent the optimality and feasibility cutting planes, respectively. It is worth noting that when the probabilities for a certain scenario are zeros in the optimality and feasibility cut planes, these two cutting planes are relaxed. Moreover, they could be transformed by strong duality or KKT   condition to solve the bilinear issues.

\subsection{Column Generation: Solving $\bf{WCEV}$ From Primal Perspective}
This reformulation in Proposition \ref{p1} naturally facilitates the application of decomposition algorithms, such as column generation (CG). By iteratively capturing the worst-case probability distribution and enlarging the corresponding scenario set, the method effectively tightens the bound while preserving computational tractability. The decomposition separates the $\mathbf{WCEV}$ problem into a pricing master problem ($\mathbf{PMP}$) and a pricing subproblem ($\mathbf{PSP}$). The $\mathbf{PMP}$ optimizes over a finite set of scenarios $\{{\bm{\xi}}_j^{\mathrm{f}}\}_{j=1}^n$, yielding the following formulation:
\begin{eqnarray}
	&&\!\mathbf{PMP}\!:\!\eta(\!\bm{\lambda},\bm{\xi}^\mathrm{f}) \! =\! \big\{\! \max_{p_{{\bm{\xi}}^{\mathrm{f}}}}\! \sum_{j=1}^n \! Q(\bm{\lambda},\bm{\xi}_j^{\mathrm{f}})p_{{\bm{\xi}_j^{\mathrm{f}}}} \!:\!\{(\bm{\bm{\xi}}^{\mathrm{f}}_j,p_{{\bm{\xi}}^{\mathrm{f}}_j})\!\}_{j\!=\!1}^n \!\in\! \mathcal P(\bm \lambda)\!\big\} \label{PMP}
\end{eqnarray}

\indent Let $\pi^*$ and $\bm{\kappa}^*$ denote the optimal dual variables associated with the constraints in \eqref{PMP}, respectively. To identify the scenario with the largest reduced cost, the $\mathbf{PSP}$ is formulated as:
\begin{eqnarray}
	&& \mathbf{PSP}:\varphi(\bm{\lambda},\pi^*,\bm{\kappa}^*) = 	\max_{\bm{\xi}\in\Omega} Q(\bm{\lambda}, \bm{\xi}) - \pi^* - \bm{\psi}(\bm\xi)^\mathrm{T}\bm{\kappa}^*
\end{eqnarray}
After linearizing $\mathbf{PSP}$ using KKT conditions, the next result, which guarantees the finite convergence of CG algorithm, follows from the properties of linear program directly.
\begin{prop}
	\label{p2}
	Suppose that  $Q(\bm{\lambda}, \bm{\xi})$ is always feasible for any $\bm{\lambda} \in \Lambda$. Given that $\Omega$  is a polytope and $\bm{\psi}$ is linear, for a fixed $\bm \lambda$, the CG algorithm terminates with an optimal solution of $\mathbf{WCEV}$ after a finite number of iterations, bounded by the number of extreme points derived in $\Omega$.
\end{prop}

\subsection{Enhancements}

\noindent $\blacksquare$ \underline{\emph{Positive Probability Scenario Filtering (PPSF)}} . 
A critical challenge in the our algorithm is the large number of scenarios generated by CG algorithm in each iteration. However, not all scenarios contribute to the worst-case distribution. To address this, we introduce PPSF, a mechanism that selectively retains scenarios with strictly positive probabilities, i.e., scenarios in the support under given $\bm \lambda$. Since the probability information can be directly obtained from $\mathbf{PMP}$, this mechanism is trivial to implement. The following result simplifies $\mathbf{MMP}$ via PPSF under certain trivial conditions.

\begin{prop}
	\label{PPSF}
	If there exists an $\epsilon>0$ such that $ \overline{\omega}_{i,j}-\varepsilon_{ij}g_{ij} \geq\underline{\omega}_{i,j} +\epsilon$ for any $\bm g$ and $(i,j)\in\mathcal{B}$, $\mathbf{MMP}$ is equivalent to the formulation below: 
	\begin{eqnarray}
		&&\textstyle \mathbf{MMP}_\mathrm{e}: \ \Gamma = \min_{\bm{\lambda}\in\Lambda} \ \bm{c}^\mathrm{T}\bm{\lambda}+\eta \label{MMP1_PPSF} \\ 
		&&\textstyle \mathrm{s.t.}\ \  \eta \geq \max\{\sum_{\bm\xi\in \Xi}(\bm{g}^\mathrm{T}\bm{\zeta}_{\bm\xi})p_{\bm{\xi}}: (p_{\bm\xi})_{\bm\xi\in\Xi}\in\mathcal{P}(\bm \lambda) \} \label{MMP2_PPSF}\\ 
	&&\textstyle  \qquad\bm\zeta_{\bm\xi}\in{\mathcal{Z}}(\bm\lambda,\bm{\xi}), \ \forall \bm{\xi}\in\Xi_\mathrm{s} \label{MMP4_PPSF}
\end{eqnarray} 
where $\Xi_\mathrm{s}$ is the scenario set generated by using PPSF.
\end{prop}

\begin{remark}
In realistic distribution systems, grid reinforcement measures cannot fully eliminate uncertainty effects; therefore, they can always be sufficiently to ensure the $\epsilon$-separation holds. Consequently, Proposition \ref{PPSF} provides a highly practical and widely applicable reformulation tool in real-world implementations.
\end{remark}

\noindent $\blacksquare$ \underline{\emph{Initialization Approach}}   
Initializing the CG procedure leads to a huge challenge due to the complexity of the ambiguity set $\mathcal{P}(\bm \lambda)$. Specifically, the ambiguity set consists of multiple wind intensity levels, whose associated bounds for the first-order moment inequalities differ from one another. 
To address this issue, we develop a specialized initialization strategy that constructs an initial feasible scenario set by leveraging extreme rays derived from the Farkas Lemma. 
Given an extreme ray ($\pi', \bm{\kappa}'$) from the dual problem of the current infeasible $\mathbf{PMP}$, the corresponding initialization subproblem is defined as:
\begin{eqnarray}
&& \mathbf{PSP}_{\mathrm{IA}}:\varphi_{\mathrm{IA}}(\bm{\lambda},\pi',\bm{\kappa}') = 	\max_{\bm{\xi}\in\Omega} - \pi' - \bm{\psi}(\bm\xi)^\mathrm{T}\bm{\kappa}'
\end{eqnarray}

This process is repeated until the updated $\mathbf{PMP}$ becomes feasible. Then, CG algorithm begins with this initialized scenario set.

\begin{algorithm}[H]
\setlength{\baselineskip}{0.8\baselineskip}
\caption{Primal-based Decomposition Algorithm (C\&CG-DRO)}
\begin{algorithmic}[1]
\STATE \textbf{Initialization:} Set $LB\!=\!-\infty,UB\! =\! +\infty, \!iter=\!1$, $\Xi_\mathrm{s} = \emptyset$, and tolerance $GAP_\epsilon$.
\WHILE{$UB - LB > GAP_\epsilon$} 
    \STATE Compute the main master problem $\mathbf{MMP}_\mathrm{e}$. 
	\IF{ $\mathbf{MMP}_\mathrm{e}$ is infeasible}
		\STATE Report the infeasibility of \textbf{DDU-DRO} and terminate.
	\ELSE
		\STATE Derive optimal solution $\hat{\bm\lambda}$, value $\Gamma^*$, and update $LB = \Gamma^* $. %
	\ENDIF
    \STATE Solve the $\mathbf{WCEV}(F)$ via CG Algorithm, derive its optimal value $\eta_\mathrm{f}(\hat{\bm\lambda})$,  and PPSF-filtered scenario set $\hat{\Xi}_{\mathrm{f}}(\hat{\bm \lambda})$
	\IF{$\eta_\mathrm{f}(\bm\hat{ \lambda}) = 0$}
    		\STATE Solve the $\mathbf{WCEV}$  via  CG Algorithm, derive its optimal value $\eta_\mathrm{o}(\hat{\bm \lambda})$,  and PPSF-filtered scenario set $\hat\Xi_{\mathrm{o}}(\hat{\bm \lambda})$. 
	\ELSE
		\STATE Set $\eta_o(\hat{\bm{\lambda}}) = +\infty$ and $ \hat\Xi_\mathrm{o}= \emptyset$.
	\ENDIF
	\STATE Update $\Xi_\mathrm{s}=\Xi_\mathrm{s}\cup\hat\Xi_{\mathrm{f}}(\hat{\bm\lambda})\cup\hat\Xi_{\mathrm{o}}(\hat{\bm\lambda})$, and augment $\bf{MMP}$ with new variables \eqref{MMP4_PPSF} and constraints \eqref{MMP2_PPSF} accordingly; 
	\STATE  Update $UB \!=\! \min\{UB, \bm{c}^\mathrm{T}\bm{\hat{\lambda}}+\eta_o(\hat{\bm{\lambda}})\}$, and $iter=iter+1$.
\ENDWHILE\\
\STATE\textbf{Output: } Report the optimal solution $\hat{\bm\lambda}$ of $\mathbf{DDU-DRO}$
\end{algorithmic}
\end{algorithm}

\section{Numerical Studies}
\label{NS}
To validate the effectiveness of the proposed $\mathbf{DDU-DRO}$ planning framework for multi-microgrids in island systems, we conduct numerical experiments on a constructed test system consisting of two load islands and one resource island. This testbed is adapted from an exemplary  distribution network and modified to emulate the operational characteristics in sea islands. 
All algorithms have been implemented in Python (Version 3.12) and executed on a Mac equipped with an Apple M4 processor (10-core CPU, 10-core GPU) running macOS 15.2 with 16 GB RAM. The master- and sub- problems are solved using Gurobi Optimizer (Version 12.0.0).

\subsection{Experiment Settings}

The test system consists of one resource island and two load islands with no inner-island cable connections. The resource island includes 6 candidate sites for the deployment of HP\&S systems. The two load islands are equipped with distribution networks containing 9 and 6 nodes, respectively.  Among these, 3 nodes on load island 1 and 1 node on load island 2 are potential sites of local renewable generation and operate FC-based H2P systems supported by  ESSs. Each load island operates independently and is equipped with power electronic interfaces and logistics infrastructure to support hydrogen reception, energy conversion, and local dispatch.

To characterize the impact of wind variability and system-level faults on long-term planning and resilience, 6 discrete levels  are defined, each representing a representative condition  based on wind speed ranges. These levels are grouped into three broad categories.
\begin{itemize}[left=0em] 
	\item  \emph{Normal Conditions}: Levels 1–3 represent typical operating conditions characterized by increasing wind resource availability, where Level 1 corresponds to partial turbine cut-in with limited generation, Level 2 to full turbine cut-in with low but stable output, and Level 3 to near-rated wind speeds enabling full-capacity operation under nominal network conditions.
	\item \emph{Degradation Conditions}: Level 4 represents a partially degraded state, where wind generation remains near-rated as in Level 3, but minor transmission faults begin to affect network integrity.
	\item \emph{Contingency Conditions}: Levels 5–6 represent low-probability and high-impact $N-k$ contingency conditions. Level 5 involves reduced wind generation, multiple line disruptions, and partial hydrogen transport disruptions. Level 6 reflects severe system failures, including extensive wind turbine shutdowns, widespread network faults, and major interruptions in hydrogen delivery.
\end{itemize}

To systematically evaluate the planning robustness and operational resilience of the proposed system, we design and compare three test cases, each representing a distinct planning strategy under uncertainty:
\begin{itemize}[left=0em] 
	\item {\bf Case 1}: Capacity planning is performed considering only normal and degradation conditions (Levels 1–4). 
	\item  {\bf Case 2}: Capacity planning is performed considering all conditions (Levels 1–6). 
	\item   {\bf Case 3}: Capacity planning is performed considering all conditions (Levels 1–6), with  pre-disaster mitigation measures including hydrogen  buffering storage and grid hardening investments.
\end{itemize}

Moreover, to consistently evaluate the resilience of each planning configuration, we assess its performance across the entire uncertainty space via post-optimization model. 
We introduce a quantitative resilience metric based on the worst-case expected value of lost load (VOLL). Given a first-stage planning decision 
, the post-evaluation is formulated as the following bilevel problem:
\begin{eqnarray}
	&&\textstyle \mathbf{WCEVOLL}: \max_{\mathbb{P}\in\mathcal{P}(\bm{\lambda})}  \int_{\bm{\xi} \in\Omega}\bar{Q}(\bm{\lambda},\bm{\xi})\mathbb{P}(d\bm{\xi}) \\
	&&\textstyle\text{where}\ \bar{Q}(\bm{\lambda}, \bm{\xi}) = \min_{\bm{\zeta}\in\tilde{\mathcal{Z}}(\bm\lambda,\bm\xi)}  \sum_{j\in\mathcal{N}}\sum_{t\in\mathcal{T}}O^{\mathrm{ls}}_{j}p_{\mathrm{ls}, j, t} 
\end{eqnarray}
where $O^{\mathrm{ls}}_{j}$ denotes the penalty coefficient associated with VOLL.

The technical specifications of key system components  are summarized in the Appendix.

\subsection{Comparative Analysis of Different Planning Configurations}
\begin{figure}[!t]
	\vspace{-10pt}
	\centering
	\includegraphics[width=4in]{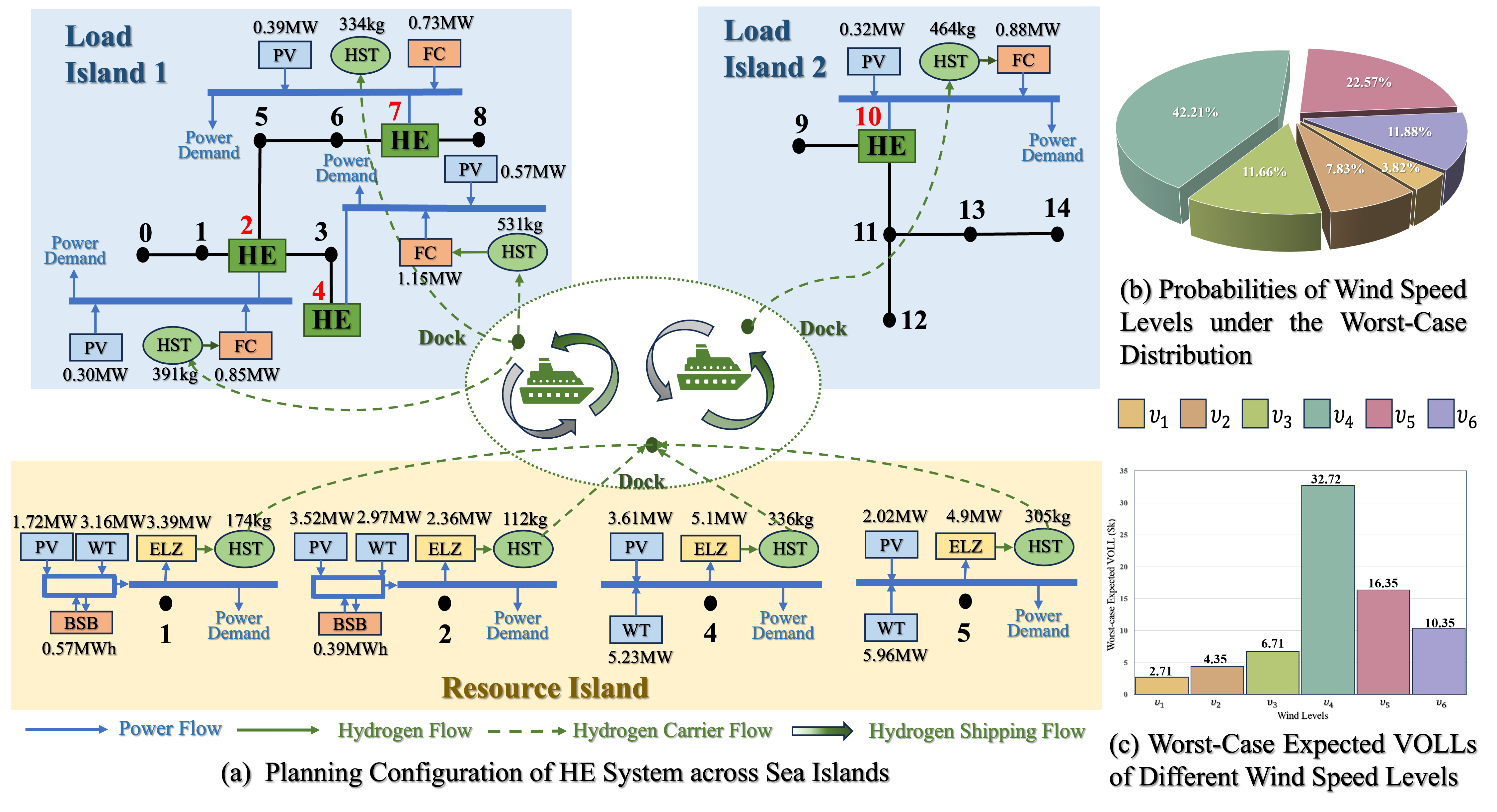}
	\caption{Planning configuration and resilience evaluation of Case 1, i.e., system deployment, hydrogen-power flows, worst-case distribution and load shedding performance} 
	\label{fig:case1} 
\end{figure}

\subsubsection{Case 1}The planning configuration of Case 1, as shown in Fig. \ref{fig:case1}(a), results in an annualized CAPEX of \$8405.27k and a worst-case expected O\&M cost of \$1940.99k. This configuration involves the deployment of four HP\&S units located at N1–N2 and N4–N5 on the resource island. 
To balance production capacity and enhance operational flexibility, these nodes are strategically organized into two spatial clusters: N4 and N5 serve as the primary hydrogen production hubs, while N1 and N2 form an auxiliary cluster that provides additional flexibility and redundancy.
Specifically, N4 and N5 are each equipped with the system's largest electrolyzers (rated at approximately 5 MW), supported by extensive RESs (e.g., 5.23MW of WTs  at N4 and 5.96MW of WTs at N5). In contrast, N1 and N2 are  equipped with medium-scale renewable generation assets and smaller electrolyzers.  For instance, N2 hosts 3.52MW of PV and 2.97MW of WTs, along with an electrolyzer rated at 3.39MW. 
The hydrogen produced on the resource island is subsequently transported to the two load islands. 
Vessel 1 serves Load Island 1, transporting hydrogen to HE nodes (N2, N4, and N7 in load island), each co-located with HSTs and FCs. 
Meanwhile, Vessel 2 services Load Island 2, delivering hydrogen directly to HE Node 10, which adopts a similar configuration for local storage and fuel cell–based power generation.
Both vessels complete two round trips per planning cycle, maintaining adequate hydrogen inventory to ensure continuous FC operation.
Additionally, small-scale PV systems installed at the HE nodes, such as 0.32 MW at N12, provide supplementary power to meet local demands. 

Then, we evaluate its resilience across the entire ambiguity set, as illustrated in Fig. \ref{fig:case1}(b) and Fig. \ref{fig:case1}(c). The total worst-case expected VOLL amounts to \$73.19k. %
A substantial portion of the expected VOLL is concentrated in Level 4, which represents scenarios with high wind power availability accompanied by minor grid faults. Such conditions are typical in real-world islanded microgrids, where operational disturbances occur frequently (though not severe). Level 4 is assigned a high worst-case probability of 0.42 and contributes \$32.72k to the total VOLL. This indicates that even light network faults under frequent operating states can lead to noticeable load shedding, making Level 4 the dominant contributor to system-level expected losses.
In contrast, Level 5 captures compound degradation scenarios, where both hydrogen vessel availability is reduced and grid damage becomes more pronounced. With a worst-case probability of 0.23, it contributes \$16.35k to the worst-case expected VOLL. Although less frequent than Level 4, Level 5 demonstrates how coupled disruptions in energy transport and infrastructure integrity can significantly elevate risk exposure, even without crossing into full-scale disaster territory.
Level 6 scenarios are characterized by severe wind generation losses, maritime hydrogen transport breakdowns, and multi-node isolation. Despite a probability of 0.12, similar to that of Level 3, Level 6 produces a substantially higher expected VOLL of \$10.35k, compared to \$6.71k for Level 3. This disparity reveals the disproportionate impact of rare but cascading failures, which can result in several times the system loss of typical disturbances, thus exposing structural vulnerabilities not fully addressed by this planning strategy.

\begin{figure}[!t]
	\vspace{-10pt}
	\centering
	\includegraphics[width=4in]{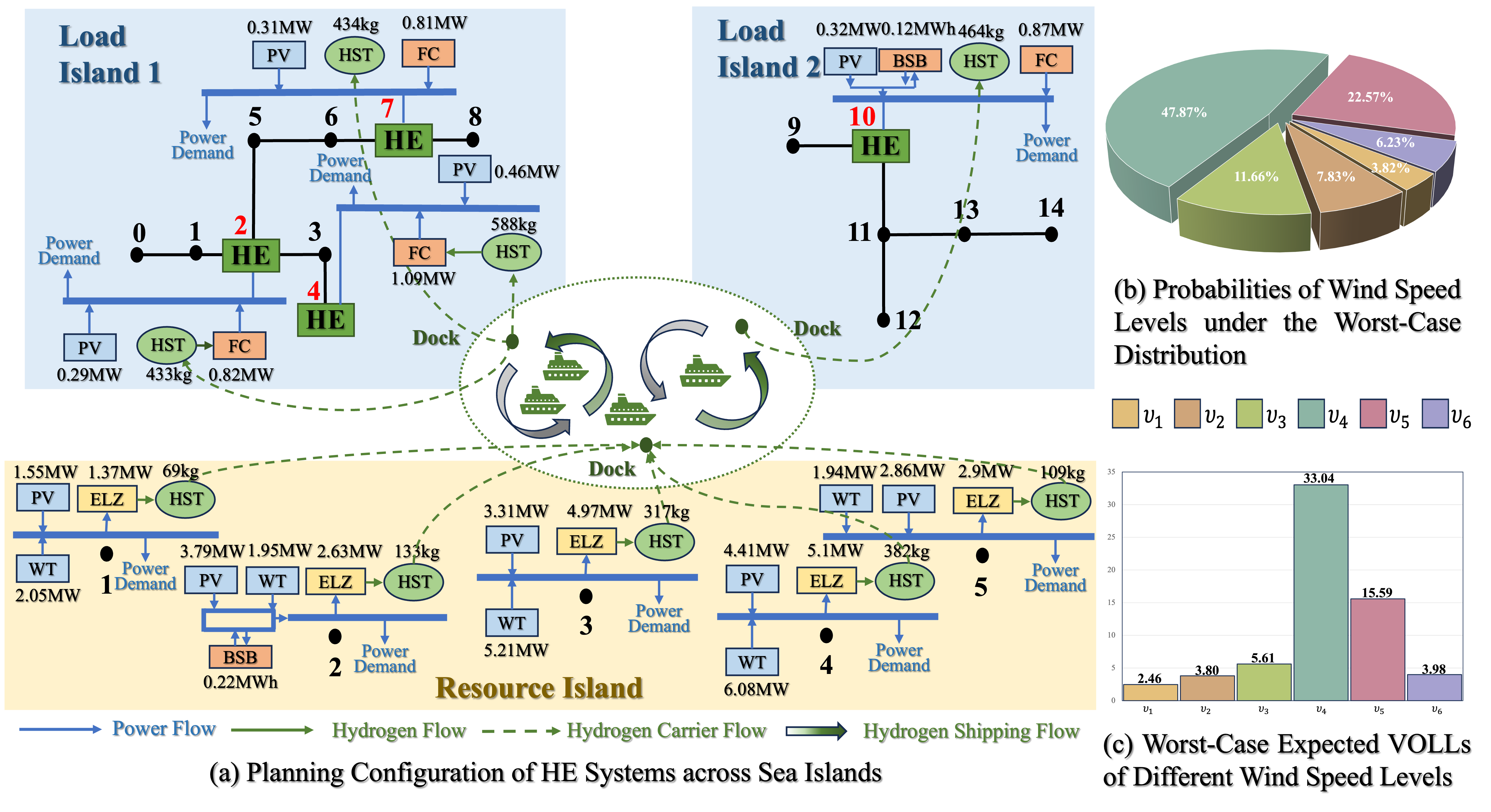}
	\caption{Planning configuration and resilience evaluation of Case 2 under full-level uncertainty} 
	\label{fig:case2} 
\end{figure}

\subsubsection{Case 2} Compared to Case 1, structural upgrades to both the production and transport subsystems are incorporated in Case 2, as shown in Fig. \ref{fig:case2}(a), with the aim of improving resilience against renewable variability and network contingencies. Specifically, an additional HP\&S unit is deployed at N3 on the resource island. This added production point is equipped with 2.9 MW of electrolyzers, 2.86 MW of PV, and 1.94 MW of WTs, with a dedicated hydrogen storage tank of 109 kg. The expanded production can reduce the risk of production bottlenecks under local renewable outages. 
Then, the hydrogen transport strategy is significantly enhanced. The number of vessels is increased from two to four, allowing for a more flexible and robust delivery network. In this configuration, two vessels are assigned to Load Island 1, one vessel to Load Island 2, and one additional vessel operates across islands, providing adaptive transport support during peak demands or supply disruptions. This flexible routing capability mitigates the effects of temporary vessel unavailability and supports load rebalancing during high-risk scenarios.
These structural improvements result in a moderately increased capital investment, with a total annualized CAPEX of \$9684.17k. However, the worst-case expected O\&M decreases to \$1842.38k, reflecting enhanced resilience and reduced penalty costs under uncertain realizations.

On the other hand, the resilience performance of Case 2 is illustrated in Fig. \ref{fig:case2}(b) and Fig. \ref{fig:case2}(c). The total worst-case expected VOLL is reduced to \$64.48k, representing a 12\% decrease compared to that in Case 1. Although the worst-case probability of Level 6 (representing compound disruptions), remains nearly unchanged, its worst-case expected VOLL decreases significantly to \$3.98k (61\% reduction). This substantial decrease verifies that the planning of Case 2 effectively mitigates cascading failures and enhances system resilience under extreme conditions. 
Meanwhile, Level 4 contributes a slightly higher worst-case expected VOLL of \$33.04k, representing an increase of only 1.0\% compared to that in Case 1.  This is associated with a moderate increase in its worst-case probability from 0.42 to 0.48. This increase reflects a trade-off made in Case 2’s planning that accepts minor efficiency losses under common conditions in exchange for significant improvements in rare but high-impact scenarios. These results demonstrate that full-level integration during the planning stage enables more resilient performance, particularly under extreme uncertainty, without substantially compromising normal-case operations.

\begin{figure}[!t]
	\vspace{-10pt}
	\centering
	\includegraphics[width=4in]{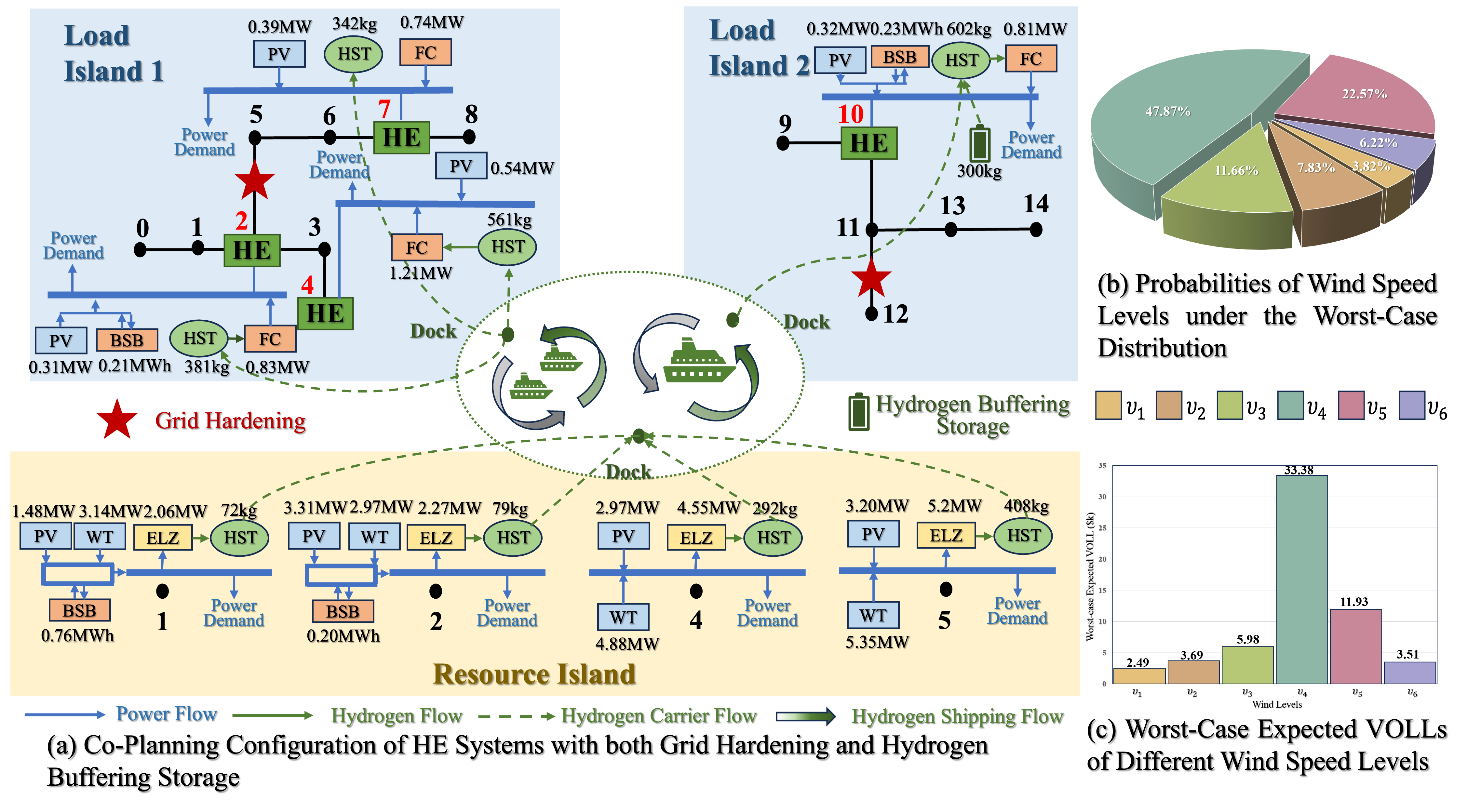}
	\caption{Planning configuration and resilience evaluation of Case 3 with grid hardening and hydrogen buffering} 
	\label{fig:case3} 
\end{figure}

\subsubsection{Case 3} Case 3 presents a planning strategy that emphasizes targeted resilience enhancements (i.e., grid protection and hydrogen buffering storage). As shown in Fig. \ref{fig:case3}(a), the  configuration  differs from Case 2 in several key aspects.   
Specifically,  the HP\&S unit  at N3 on the resource island is removed, reducing the number of active production nodes from five to four. Also, the number of hydrogen vessels is downsized from four to three, with two vessels assigned to Load Island 1 and one to Load Island 2.  These adjustments streamline the infrastructure, resulting in a reduction of the annualized CAPEX to \$8382.27k, representing a 13.4\% decrease compared to Case 2.  
While infrastructure scale is reduced, two structural resilience measures are introduced to enhance system robustness. Specifically, grid hardening is applied to Line 2-5 and Line 11-12. Second, a hydrogen buffering storage system is deployed at N10 on Load Island 2, providing local supply stability in the event of vessel delays or upstream disruptions.
As a result, these targeted enhancements contribute to a reduction in the worst-case expected O\&M cost, which drops to \$1711.14k despite the streamlined system configuration.

As illustrated in Fig. \ref{fig:case3}(b) and Fig. \ref{fig:case3}(c), the probability distribution of load shedding  in Case 3 remains largely similar to that in Case 2, indicating no significant degradation in overall system resilience.
However, the total worst-case expected VOLL is further reduced to \$60.98k, representing a 5.4\% decrease compared to Case 2. This improvement reflects the effectiveness of the grid hardening and hydrogen buffering strategies in mitigating local high-impact disruptions, particularly under Level 5 and Level 6 scenarios. By reinforcing critical transmission pathways and introducing local hydrogen reserves, Case 3 enhances the system’s resilience against compound failures without increasing risk under more probable operating conditions.

In summary, Case 3 demonstrates that strategic infrastructure streamlining, coupled with selective resilience enhancements, can achieve both cost savings and improved system performance. By reducing production and transport capacity while reinforcing critical transmission pathways and introducing local hydrogen reserves, the system attains a more effective balance among investment cost, operational efficiency, and resilience.

\subsection{Algorithm performance}
To evaluate the computational effectiveness of our proposed solution approaches, we compare two algorithms: the basic-C\&CG algorithm with strong cutting planes and the  C\&CG-DRO algorithm, which incorporates both strong cutting planes and enhancements within the CG decomposition framework. 

The comparison is performed on three synthetic systems consisting of 5, 10, and 15 buses, respectively. For each system size, we simulate three levels of load intensity: 0.5, 0.7, and 0.9 (normalized). The results are summarized in Table \ref{AP}. For each algorithm and test setting, we report the lower bound (LB), upper bound (UB), relative optimality gap, number of iterations (Iter.), number of distinct scenarios generated (Scen.), and total computation time in seconds. All experiments are conducted under identical initial conditions to ensure a fair comparison. The computational results reveal a consistent performance advantage of the C\&CG-DRO algorithm over basic-C\&CG across all system sizes and load levels. 
\begin{table*}[!t]
	\centering
	\footnotesize  %
	\caption{Comparison Between basic-C\&CG and C\&CG-DRO Under Different Bus Numbers and Load Levels}
	\renewcommand{\arraystretch}{0.9}
	\setlength{\tabcolsep}{3pt}
	\newcolumntype{Y}{>{\centering\arraybackslash}X}  
	
	\begin{tabularx}{0.9\textwidth}{YY|YYYYYY|YYYYYY}
		\toprule
		\multicolumn{2}{c|}{\textbf{System Settings}} & \multicolumn{6}{c|}{\textbf{basic-C\&CG}} & \multicolumn{6}{c}{\textbf{C\&CG-DRO}} \\
		\cmidrule(lr){1-2} \cmidrule(lr){3-8} \cmidrule(lr){9-14}
		\textbf{Bus No.} & \textbf{Load} & \textbf{LB} & \textbf{UB} & \textbf{Gap} & \textbf{Iter.} & \textbf{Scen.} & \textbf{Time (s)} & \textbf{LB} & \textbf{UB} & \textbf{Gap} & \textbf{Iter.} & \textbf{Scen.} & \textbf{Time (s)} \\
		\midrule
		5  & 0.5 & 3816.38 & 3819.43 & 0.08\% & 46 & 46 & 1554.10  & 3817.32    & 3818.08 & 0.02\% & 5 & 38 & 221.07 \\
		& 0.7 & 4202.29 & 4205.23 & 0.07\% & 46 & 46 & 2139.15 & 4203.85 & 4204.27 & 0.01\% & 7 & 43 & 292.24 \\
		& 0.9 & 4568.56 & 4570.84 & 0.05\% & 53 & 53 & 2582.89 & 4569.28 & 4571.56  & 0.05\% & 8 & 54 & 350.06 \\
		\midrule
		10 & 0.5 & 5005.25 & 5008.75 & 0.07\% & 67 & 67 & 8007.81 & 5005.58 & 5007.58 & 0.04\% & 6 & 59 & 914.64 \\
		& 0.7 & 6241.27 & 6246.88 & 0.09\% & 65 & 65 & 6150.99 & 6243.07 & 6244.31  & 0.02\% & 5 & 58 & 944.18 \\
		& 0.9 & 7485.46 & 7492.19 & 0.09\% & 66 & 66 & 9642.37 & 7488.95 & 7491.19  & 0.03\% & 5 & 60 & 843.31 \\
		\midrule
		15 & 0.5 & 7546.39 & 7553.18 & 0.09\% & 59 & 59 &12565.28 & 7545.06 & 7548.83  & 0.05\% & 7 & 47 & 1231.39 \\
		& 0.7 & 9561.81& 9569.45 & 0.08\% & 50 & 50 & 8104.53 & 9561.03 & 9562.94  & 0.02\% & 7 & 55 & 984.91 \\
		& 0.9 &12874.38 &12882.10 & 0.06\% & 53 & 53 &11926.35 &12874.50 &12879.65  & 0.04\% & 6 & 40 & 612.23 \\
		\bottomrule
	\end{tabularx}
	\label{AP}
\end{table*}

1) \emph{Solvability and Stability}:
Both algorithms successfully solved all instances, demonstrating the baseline tractability of the proposed methods. However, C\&CG-DRO shows significantly more stable convergence behavior, with convergence typically achieved within 5–8 iterations, compared to  iterations (up to 67) for basic-C\&CG. This indicates that C\&CG-DRO better captures the problem structure and avoids redundant scenario exploration. 

2) \emph{Solution Quality}: 
The optimality gap achieved by C\&CG-DRO is consistently lower than that of basic-C\&CG. In smaller systems (e.g., 5-bus), the gap is reduced from 0.08\% to 0.02\% under low load; in larger systems (15-bus), gap improvements remain notable (e.g., 0.09\% to 0.04\%). This demonstrates the higher precision and robustness of the proposed method in approximating the worst-case distribution, even under high-stress operational conditions.

3) \emph{Algorithmic Efficiency}:
In terms of scenario economy, C\&CG-DRO typically requires 20–30\% fewer scenarios than basic-C\&CG to achieve tighter bounds. This reduction directly translates to a drastic drop in computation time. For instance, in the 15-bus system at 0.9 load, C\&CG-DRO reduces runtime from 11,926s to 612s (a  95\% reduction) while producing comparable or better solutions. The embeded oracles for solving $\mathbf{WCEV}$ problem efficiently filters out non-informative scenarios early, focusing computation on critical worst-case regions.

4) \emph{Scalability Under Increasing System Complexity}:
As the system size and load level increase, the computational burden of basic-C\&CG scales sharply, with runtime and iterations rising nonlinearly. In contrast, C\&CG-DRO maintains relatively stable performance with mild growth in scenario count and time, reflecting its better scalability and structural adaptivity. This becomes especially valuable for large-scale planning problems under complex uncertainty.

In summary, while both algorithms are capable of solving the problem, C\&CG-DRO demonstrates superior performance across all critical dimensions. It achieves tighter solution bounds with fewer scenarios, scales more favorably with system complexity, and significantly reduces computation time by an order of magnitude on average. These results validate the effectiveness of the primal-based decomposition framework in enabling efficient and resilient infrastructure planning under distributional uncertainty.

\section{Conclusion}
\label{CONC}
This paper proposes a two-stage distributionally robust planning framework for HE sea-island multi-microgrid systems. A maritime hydrogen transport network is integrated to enable spatial-temporal decoupling between renewable power generation  and consumption across islands. Comprehensive DDU-based wind-driven uncertainties are incorporated into the model. A primal-based decomposition algorithm with strong cutting planes is developed to  efficiently solve the DDU-DRO model, without a complete recourse assumption.
\begin{enumerate}
	\item \emph{Trade-Off between Cost Efficiency and System Resilience in Strategic Planning}:
	Planning decisions solely considering normal operational scenarios are insufficient to ensure resilience under deep uncertainty. Organizations must adopt scenario-diverse planning frameworks that explicitly consider rare but high-impact disruptions during the decision-making stage. By embedding broad uncertainty modeling into early planning, managers can mitigate the risks of catastrophic failures without incurring excessive operational penalties during normal periods.
	\item \emph{Selective Reinforcement over Uniform Infrastructure Expansion}: Resource allocation strategies that target critical vulnerabilities, such as reinforcing key network lines or providing local backup capabilities, offer superior resilience-to-cost ratios compared to uniform infrastructure expansion. The findings emphasize that in the planning of isolated sea-island  microgrids, decision-makers should prioritize selective, risk-informed investments rather than pursuing undifferentiated system growth, thus achieving more sustainable and robust operational outcomes.
	\item \emph{Efficiency and Scalability of the Proposed Algorithms}: 
	The computational results demonstrate that the basic-C\&CG algorithm with strong cutting planes solves the distributionally robust planning problem under uncertainty. Nevertheless, C\&CG-DRO achieves a great improvements in convergence speed and optimality gap reduction, thereby significantly accelerating the solution process and enhancing practical scalability for large-scale system applications.
\end{enumerate}
\bibliographystyle{IEEEtran}
\bibliography{bibliography}

\ifCLASSOPTIONcaptionsoff
  \newpage
\fi

\end{document}